\documentclass[reqno,11pt]{amsart}
\usepackage{amsmath,amsthm,amscd,amsfonts,amssymb,hyperref,floatrow}
\usepackage{graphicx, color,dsfont,xcolor}

\numberwithin{equation}{section}

\newtheorem{theorem}{Theorem}[section]
\newtheorem{lemma}[theorem]{Lemma}
\newtheorem{proposition}[theorem]{Proposition}

\newtheorem{rem}[theorem]{Remark}

\newtheorem{conjecture}[theorem]{Conjecture}

\DeclareMathSymbol{\leqslant}{\mathalpha}{AMSa}{"36} 
\DeclareMathSymbol{\geqslant}{\mathalpha}{AMSa}{"3E} 
\DeclareMathSymbol{\eset}{\mathalpha}{AMSb}{"3F}     
\renewcommand{\leq}{\;\leqslant\;}                   
\renewcommand{\geq}{\;\geqslant\;}                   

\def \C{ \mathbb  C }
\def \H{ \mathbb  H }
\newcommand{\R}{\mathbb{R}}

\newcommand{\N}{\mathbb{N}}

\def \P{ \mathbb P  }

\begin{document}
\title{Random matrices in non-confining potentials}

\author{Romain Allez \and Laure Dumaz}
\address{Weierstrass Institute, Mohrenstr. 39, 10117 Berlin, Germany.}
\address{Statistical Laboratory, Centre for Mathematical Sciences, Wilberforce Road, Cambridge, CB3 0WB, United Kingdom. }

\email{romain.allez@gmail.com, L.Dumaz@statslab.cam.ac.uk} 
\date{\today}

\maketitle

\begin{abstract}
We consider invariant matrix processes diffusing in non-confining cubic potentials of the form $V_a(x)= x^3/3 - a x, a\in \R$. 
We construct the trajectories of such processes 
for all time by restarting them whenever an explosion occurs, from a new (well chosen) initial condition, insuring 
continuity of the eigenvectors and of the
non exploding eigenvalues. We characterize the dynamics of the spectrum in the limit of large dimension and 
analyze the stationary state of this evolution explicitly. 
We exhibit a sharp phase transition for the limiting spectral density $\rho_a$
at a critical value $a=a^*$. 
If $a\geq a^*$, then the potential $V_a$ presents a well near $x=\sqrt{a}$ deep enough to confine all the particles inside, 
and the spectral density $\rho_a$ is 
supported on a compact interval. 
If $a<a^*$ however, the steady state is in fact dynamical with a macroscopic stationary flux of particles 
flowing across the system. In this regime,  
the eigenvalues allocate according to a stationary density profile  
$\rho_{a}$ with full support in $\R$, flanked with heavy tails such that $\rho_{a}(x)\sim C_a /x^2$ as 
$x\to \pm \infty$. Our method applies to other non-confining potentials and we further investigate 
a family of quartic potentials, which were already studied in \cite{bipz} to count planar diagrams. 
\end{abstract}

\section{Introduction}

Since Wigner's initial suggestion \cite{wigner2} 
that the statistical properties of the eigenvalues of random Hermitian matrices should provide 
a good description of the excited states of complex nuclei, random matrix theory (RMT) has 
become one of the prominent fields of research, at the boundary between atomic physics, 
solid state physics, statistical mechanics, statistics, probability theory, and number theory \cite{handbook}.

The two main models for Hermitian random matrices (which have been extensively studied in the literature, 
see \cite{handbook,agz,mehta,silverstein,forrester} for a review of RMT) are: 
(i) the ensembles of Wigner matrices which have independent (up to symmetry) and identically distributed entries;
(ii) the classical real, complex and quaternion {\it invariant} ensembles, which are respectively stable   
under the conjugation by the orthogonal, unitary, and symplectic groups. 
The intersection between those two types of random matrices is actually reduced to the famous Gaussian orthogonal, unitary and symplectic ensembles.

In this paper, we construct Hermitian matrix diffusion processes $(H(t),\;t\geq 0)$, which 
evolve in {\it non-confining} potentials and which are invariant under rotation at all time $t\geq 0$.  
The initial motivation for our work is to make sense of invariant ensembles of random matrices 
in non-confining potentials such as $V_a(x)= x^3/3-a\,x,\; a\in \R$, which cannot be defined in the usual way because of the divergence of the partition function.  
The family $\{V_a,a\in \R\}$ of potentials we are looking at includes the cubic interaction ($a=0$) and our results  
can be extended to a variety of non-confining polynomial potentials with the same ideas.   
In particular, our method also permits us to construct a family of invariant ensembles in non-confining 
quartic potentials of the form $U_g(x)= x^2/2+ g\, x^4$ where $g<0$ is a  parameter called the coupling constant. 
Such potentials $U_g$ have already been considered in the paper \cite{bipz}, where 
Br\'ezin, Itzykson, Parisi and Zuber use an enumerative formula established in \cite{hooft} 
for planar diagrams in terms of matrix integrals 
associated to the potential $U_g$,  to count planar diagrams.  Our construction brings new lights on some of their results 
on the limiting spectral density as the dimension $N\to \infty$, in the case  $g<0$.

The idea of our construction is inspired by \cite{halperin} where Halperin studies one dimensional diffusion 
processes in non-confining potentials. We first simply let our diffusive matrix process 
evolve in the non-confining potential of interest, 
until the first explosion time of one of the coefficients of the matrix. 
To extend the trajectory after this explosion time, we restart our matrix process at this time
from a new position insuring continuity of the non-exploding coefficients, 
until the next explosion occurs, and so on. 
This procedure is explained in section \ref{cubic-hermitian-process}.  
After some time, we expect that this matrix process will reach an equilibrium, which is, in contrast with the classical 
case of Dyson Brownian motion in a quadratic potential, a {\it dynamical steady state}. 
In particular, we establish in Section \ref{stationary-flux} that
there is a stationary flux of particles (eigenvalues) flowing across $\R$ in the steady state 
if the barrier of the potential is not too strong to confine the particles. 
This feature appears to be new and leads to interesting further questions 
about the fluctuations of this flux around its leading order in the large $N$-limit. 
We compute the leading behavior of the stationary 
flux in the large $N$ limit (see formula \ref{leading-flux}), 
which measures the number of eigenvalues crossing over the system per unit of time.  
We believe that our construction defines an interesting model of interacting charged particles 
with a random flux, to be related with the 
asymmetric exclusion processes which have attracted attention in the last decade (see e.g. \cite{ferrari,spohn}).     

In the particular case of a cubic potential, we describe precisely the dynamics of the spectrum of our matrix process in the limit of large dimension $N$ 
in section \ref{dynamic-density}. We analyze the stationary state of this dynamic by computing explicitly 
the spectral density, which gives the global allocation of the eigenvalues. We observe an interesting sharp phase transition 
for this spectral density at the critical value $a=a^*$. If $a<a^*$, then the well of the potential $V_a$  is not 
deep enough to prevent a macroscopic proportion of the eigenvalues to explode and we find a stationary 
spectral density $\rho_{a}$ with full support in $\R$, flanked with heavy tails such that $\rho_{a}(x)\sim C_a /x^2$ as 
$x\to \pm \infty$. On the contrary, if $a \geq a^*$, then the eigenvalues are confined in the well near $x=+\sqrt{a}$
of the potential $V_a$ such that the spectral density in the stationary state has compact support.
The underlying flux of particles displays also a phase transition at the critical value $a=a^*$. The current of particles, 
which allocate according to the stationary profile $\rho_a$, is 
macroscopic compared to the total mass of the $N$ eigenvalues  if $a< a^*$ whereas it is microscopic if $a\geq a^*$.

We conclude in the last section \ref{conclusion} with some open questions related to the statistics of the
eigenvalues and their current, in the stationary state.

{\bf Acknowledgments.}
We are grateful to Jo\"el Bun and Antoine Dahlqvist for interesting discussions on Stieltjes transforms. 
R.A. received funding from the European Research Council under the European
Union's Seventh Framework Programme (FP7/2007-2013) / ERC grant agreement nr. 258237 and thanks the Statslab in DPMMS,
Cambridge for its hospitality.
The work of L.D. was supported by the 
Engineering and Physical Sciences
Research Council under grant EP/103372X/1 and L.D. thanks the hospitality of the maths department of TU and the Weierstrass institute in Berlin.

\section{Invariant ensembles}
Given an {\it analytic} potential $V$ such that $V(\R)\subseteq\R$, we associate an invariant ensemble of random matrices in the space of 
(real, complex or quaternion) Hermitian matrices 
by specifying their law as
\begin{align}\label{inv-ens}
P(dH) := \frac{1}{Z}\, \exp\left(- N \, {\mbox{Tr}}(V(H))\right) \, dH\,,
\end{align}
where $V$ is meant to act on Hermitian matrices by holomorphic functional calculus, 
$dH$ is the Lebesgue measure on the space of 
Hermitian matrices and 
$Z$ is the partition function given by
\begin{align}\label{partition-function}
Z:= \int  \exp(- N \, {\textrm{Tr}}(V(H))) \, dH\,, 
\end{align}
where the integral is over the corresponding space of Hermitian matrices. 

An assumption on the growth of $V(x)$ when $x\to \pm \infty$ is necessary here \cite{biane2}, such that
\begin{equation}\label{restriction}
Z < \infty\,.
\end{equation} 
so that the partition function $Z$ and the probability law $P$ are well defined.

Under this assumption,  the eigenvalues of the matrix $H(t)$, which evolves in the stationary potential 
$V$, are confined and do not explode.

 The purpose of this article is to construct new invariant ensembles 
of random matrices in a family of potentials $(V_a, \,a\in \R)$, which violates the key assumption \eqref{restriction}. 
For simplicity, we reduce to the case of real symmetric matrices, although this discussion can be extended to the complex and quaternion cases.

The first step is to view the probability distribution $P$ defined in \eqref{inv-ens} as the {\it Boltzmann equilibrium weight} of the 
diffusive matrix process satisfying the Langevin equation  
\begin{align}\label{langevin-matrix}
dH(t) = - \frac{1}{2} \, V'(H(t)) dt + \frac{1}{\sqrt{N}} dB(t)\,,
\end{align}
with $H(0)=0$ at the initial time and
where $(B(t))$ is a {\it Hermitian Brownian motion}, i.e. a real symmetric $N\times N$ matrix process with entries given by independent 
(up to symmetry) real Brownian motions (BMs), with an extra factor $1/\sqrt{2}$ for the off diagonal entries of $(B(t))$ 
\footnote{To insure invariance by rotation, the variance of the real BMs on the diagonal has to be 
twice the off diagonal terms.}.
The dynamic \eqref{langevin-matrix} preserves the invariance under rotation of the process $H$, 
in the sense that the matrix $H(t)$ is invariant in law under the conjugation of any orthonormal matrix: 
For all $O\in \mathcal{O}_N$ (the orthonormal group) and for any $t\geq 0$, we have
\begin{align*}
H(t)  \stackrel{(d)}{=} O \, H(t) \, O^\dagger
\end{align*}
where $O^\dagger$ refers to the transpose of $O$.

 The interesting feature of this second approach is that one can still define a Hermitian matrix process
satisfying the Langevin equation \eqref{langevin-matrix} even if the potential $V$ does {\it not} fulfill 
the {\it confining} restriction \eqref{restriction}.

\section{Hermitian matrix process in cubic potentials}\label{cubic-hermitian-process}

In this paragraph, we consider the {\it odd} potential 
\begin{align}\label{cubic-potential}
V_a(x) = \frac{x^3}{3} - a \,x\,,
\end{align}
where $a\in \R$ is a given parameter. This family of potential $(V_a,\;a\in \R)$ 
includes the cubic interaction (case $a=0$) 
and we shall soon see that it is also natural to introduce a linear term in order to cover a wider variety of behaviours. 
Besides, this family of potentials already appeared in various contexts: random matrices 
\cite{virag,laure1,bloemendal1,bloemendal2}, random Schr\"odinger operators and diffusions \cite{halperin,lloyd,mckean,texier,laure2,ito}. 

\begin{figure}[h!btp] 
     \center
     \includegraphics[scale=0.9]{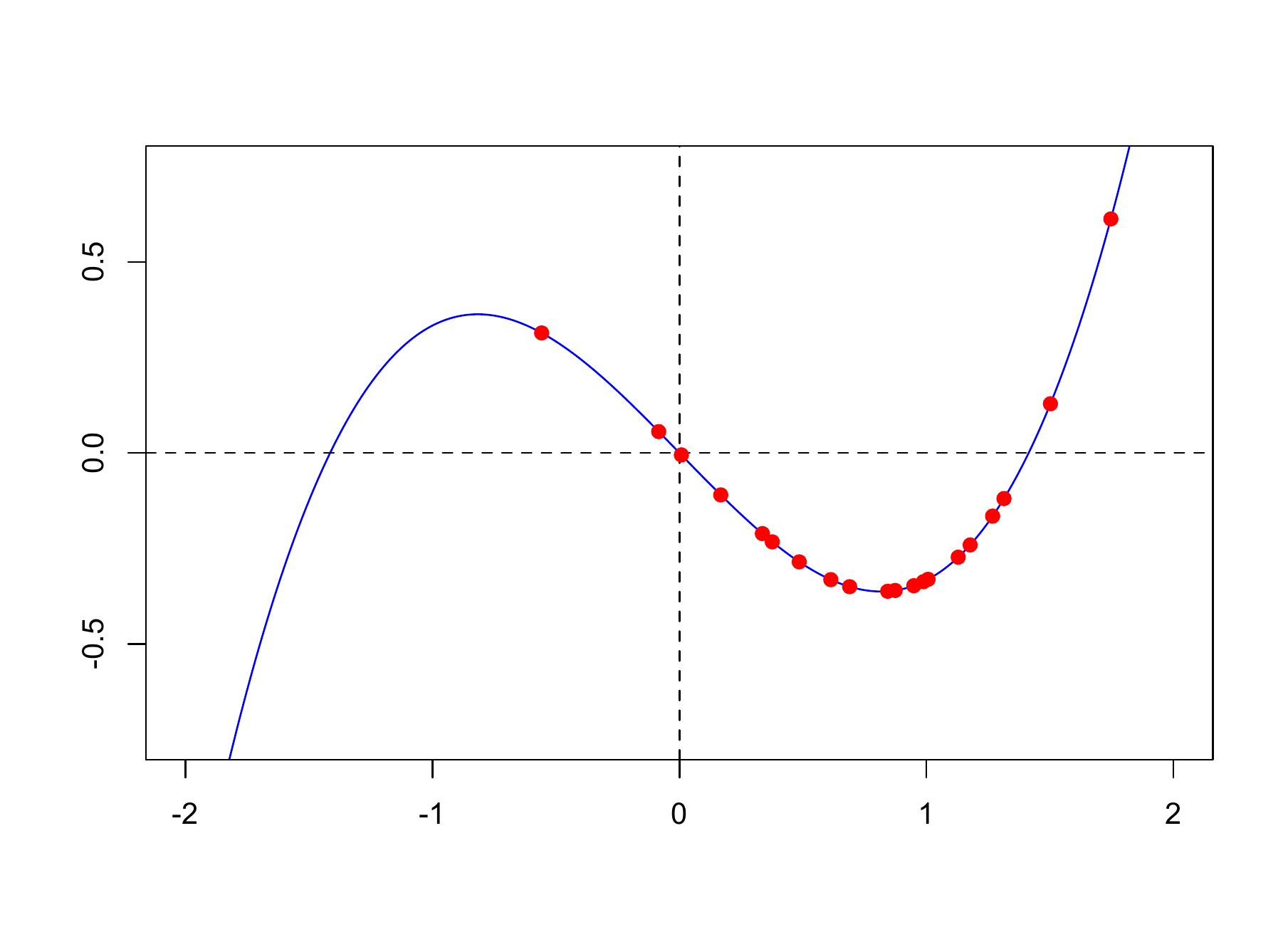}
     \caption{Picture of a Coulomb gas with $N=20$ particles (red dots) carrying positive charges in the 
     potential $V_a(x)$ as a function of $x$ for $a=2/3$.}\label{fig.potential}
\end{figure}

If $a > 0$, the potential $V_a$ presents a local minimum in 
$x=\sqrt{a}$ and a local maximum in $x=-\sqrt{a}$ (see Fig. \ref{fig.potential}). The local minimum is then separated from the local maximum by 
a potential barrier of size $\Delta V_a= \frac{4}{3} a^{3/2}$.  If $a\leq 0$ however, the potential is fully non-confining.

For such a potential $V_a$, the probability $P$ defined in \eqref{inv-ens} 
does not exist because of the divergence 
of $x^3$ as $x\to -\infty$ which prevents the partition function $Z$ to be finite. 
Nevertheless, the Langevin equation \eqref{langevin-matrix} remains well defined up to a small time interval $[0;\varepsilon]$ (with overwhelming probability) and there exists a diffusive matrix process $(H(t),\,t\in [0;\varepsilon])$ in the space of real Hermitian matrices such that, 
\begin{align}\label{langevin-cubic}
dH(t) = (a-H(t)^2) \, dt + \frac{1}{\sqrt{N}} dB(t)\,,
\end{align}
with $H(0)=0$ at the initial time and 
where $B$ is a Hermitian Brownian motion (defined above). Similar diffusion processes have also been thoroughly 
considered in dimension $N=1$ in 
\cite{halperin,virag, laure1,bloemendal1, lloyd,mckean,texier,laure2}. 
We invite the reader to look at \cite{laure2} for a brief review on the dynamics 
of a one-dimensional diffusion in such a potential.

As in the one-dimensional case, the main difficulty to define such a Hermitian process \eqref{langevin-cubic} on the whole positive half line $t\geq 0$, 
comes from the fact that, with probability one, the diffusive matrix process $(H(t))$ 
satisfying \eqref{langevin-cubic} will eventually blow-up at some finite time $\tau_1$ 
defined as 
\begin{align*}
\tau_1 := \inf\{t\geq 0: \max(|H_{ij}(t)|, 1\leq i,j \leq N ) = +\infty\}\,.
\end{align*}

In \cite{bloemendal2}, Bloemendal and Vir\'ag construct an Hermitian diffusive process $(H(t),\;t\geq 0)$ 
on the whole positive half line $t\geq 0$, satisfying a similar equation to \eqref{langevin-cubic} off what they call the {\it focal} points 
(see \cite[Eq. (5.8)]{bloemendal2}), which corresponds to the explosion times of the process $(H(t))$. 
They are interested in the case of a {\it non stationary} potentials similar to ours with an additional linear term $ r t$ in the drift ($r > 0$). Our case would simply corresponds to the $r=0$ case.  
The authors of \cite{bloemendal2} use a matrix generalization of Sturm oscillation theory, 
which goes back to the work of Morse \cite{morse1} (see also \cite{morse2,reid1,reid2,baur}).

Let us now explain how the eigenvalues and eigenvectors of $H(t)$ evolve until the first explosion time $\tau_1$ 
and see how the trajectories of those processes can be extended after this time.   
  
For any $t\geq 0$, the real eigenvalues of the Hermitian matrix $H(t)$ will simply be denoted, in non increasing order, as    
$\Lambda_N(t):=(\lambda_1(t)\geq \lambda_2(t) \geq \dots \geq\lambda_N(t))$.  
The main point is that the symmetric matrix $dB(t)$ in \eqref{langevin-cubic}
is invariant under conjugation by an orthogonal matrix so that the usual derivation of {\it Dyson's Brownian motion} \cite{dyson} 
is easily extended to this case.  
Indeed, the authors of \cite{bloemendal2} derive 
 the stochastic differential system satisfied by the eigenvalues process using Hadamart's variation formula, 
 see \cite[Theorem 5.4]{bloemendal2}, which is somehow the rigorous way  
of performing basic perturbation theory in the eigenvalues problem associated to the Langevin equation \eqref{langevin-cubic}  
\footnote{ An indirect derivation, which takes the stochastic differential system of the eigenvalues 
as granted in order to recover {\it a posteriori} the matrix equation \eqref{langevin-cubic}
can also be done as in \cite{agz} for the usual Dyson Brownian motion.}. Either way, we eventually obtain the 
following stochastic differential system for the eigenvalues $\Lambda_N(t)$,    
\begin{align}\label{sde.ev}
d\lambda_i = (a- \lambda_i^2) \, dt + \frac{\beta}{2N} \sum_{j\neq i} \frac{dt}{\lambda_i-\lambda_j}  + \frac{1}{\sqrt{N}} dB_i\,,
\end{align}
where $\beta=1$ and $B_i,i=1,\dots,N$ are real independent Brownian motions.
Note that the cases $\beta=2, 4$ may have also been covered using complex and quaternion Hermitian Brownian motions. 
Let us simply mention that for $\beta \geq 1$, the electrostatic repulsion is strong enough 
to prevent any collision between the eigenvalues so that the stochastic differential system has a well defined 
and continuous solution in the It\^o's sense \cite{agz}. 
Towards a physical picture, we can see the process $(\lambda_1,\dots,\lambda_N)$ as a one dimensional 
repulsive
Coulomb gas of $N$ positively charged particles, subject to a thermal noise and lying in the {\it non-confining cubic} potential
$V_a$ \eqref{cubic-potential} (see Fig. \ref{fig.potential}).  

\begin{figure}[h!btp] 
    \includegraphics[scale=0.9]{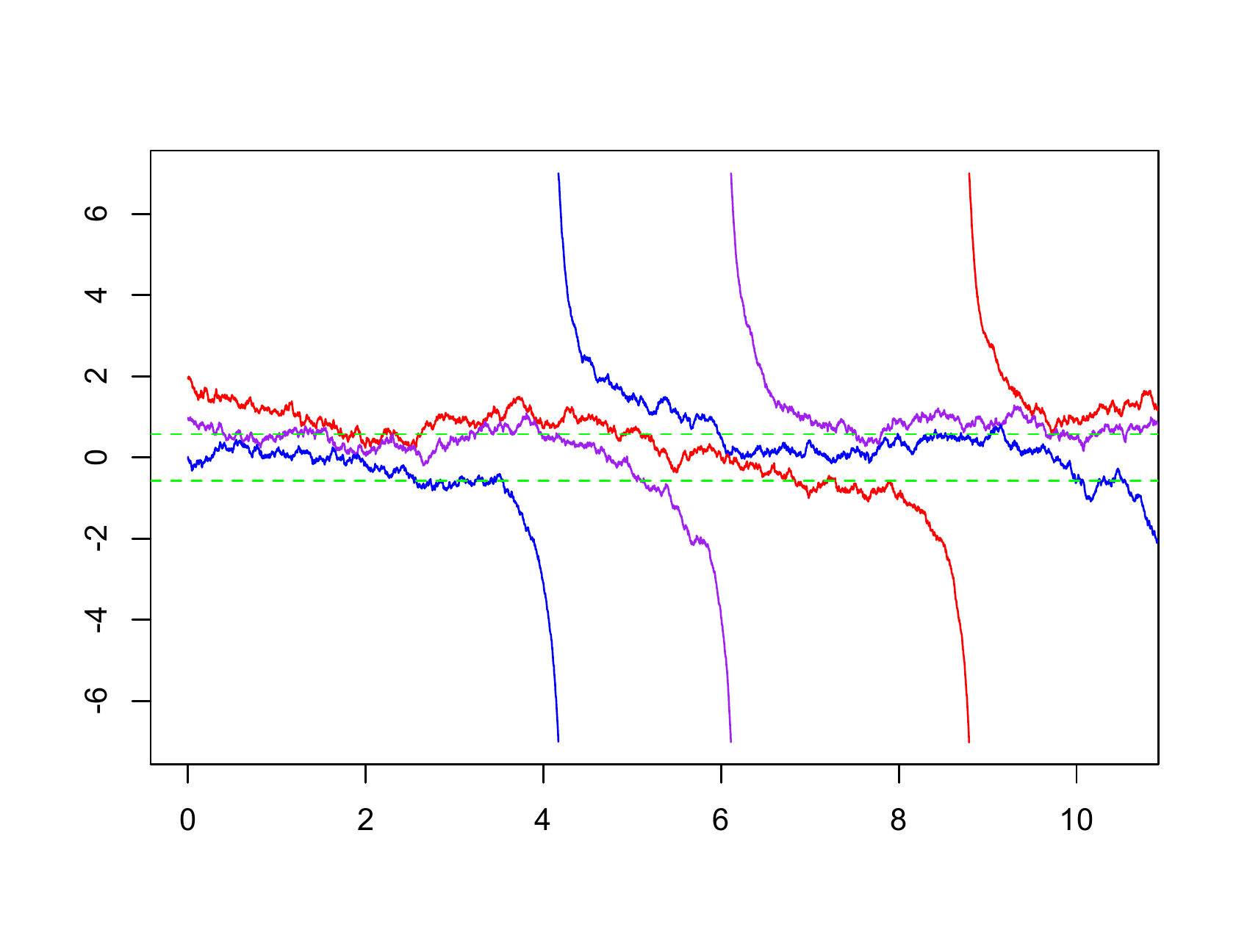}
     \caption{(Color online). Simulated paths of the eigenvalues $(\lambda_1(t),\lambda_2(t),\lambda_3(t))$ 
      as a function of time $t$ 
     for $N=3$ and $a=1/3$. The green horizontal dashed lines give the position of the well and the hill  of the potential $V_a$ in 
     $\lambda=\pm \sqrt{a}$. }\label{fig.eigenvalues-path}
\end{figure}

The evolution of the (orthonormal) eigenvectors $\psi_i(t)$ respectively associated to $\lambda_i(t)$ can also be derived 
using standard perturbation theory, or by the indirect method of \cite[Proof of Theorem 4.3.2]{agz}, 
applied for the Hermitian Brownian motion, together with the eigenvalues. 
Note that the $\psi_i(t)$ are all determined up to a sign $\pm 1$. Up to an arbitrary choice at the initial time,  we can prove, following \cite[Proof of Theorem 4.3.2]{agz} (see also \cite{alice}), 
that there exists a continuous (with respect to time) version of the process $\Psi:=(\psi_1,\dots,\psi_N)$ 
which evolves according to 
\begin{align}\label{sde-vectors}
d\psi_i(t) = - \frac{1}{2N} \sum_{j \neq i} \frac{dt}{(\lambda_i-\lambda_j)^2}\,  \psi_i(t)
+ \frac{1}{\sqrt{N}} \sum_{j \neq i} \frac{dW_{ij}(t)}{\lambda_i-\lambda_j} \, \psi_j(t)\,,
\end{align}
where the real Brownian motions $W_{ij}, 1 \leq i < j \leq N$ are mutually independent and defined by symmetry  
$W_{ij} = W_{ji}$ for $i > j$. Moreover the $W_{ij}, i\leq j$ are independent of the Brownian motions $B_i$
driving the stochastic differential system of the eigenvalues \eqref{sde.ev}.
This allows us to freeze the trajectories of the eigenvalues until the first explosion time $\tau_1$ 
and then to study the eigenvectors dynamics with this realization of the eigenvalues path. 

Now the main problem is to understand the behavior of the eigenvectors when we approach the explosion time $\tau_1$ 
at which $\lambda_N(t)\to -\infty$ as $t\to \tau_1$. 
We can easily see that this singularity does not affect the eigenvectors in the sense that, for all $i=1,\dots,N$, 
the trajectory of $\psi_i(t)$ can be extended continuously when $t\uparrow \tau_1$. 
Indeed, all the terms of the form $1/(\lambda_N-\lambda_j), j \neq N$ which appear in \eqref{sde-vectors} 
vanish at the explosion time $\tau_1$, so that we can check the Cauchy criterion 
$\psi_i(t)-\psi_i(s)\to 0$ when $t,s\to \tau_1$, which insures the existence of a limit for $\psi_i(t)$ when $t\uparrow \tau_1$ 
for all $i$.

This remark suggests to extend the trajectory of the matrix process $(H(t))$ after each explosion times, which are labeled as $\tau_k,
k\geq 1$, according to the following procedure.   
Whenever an explosion occurs at some time $\tau_k$, 
the (exploding) eigenvalue $\lambda_N$ is immediately restarted 
at the explosion time $\tau_{k+}$ 
in $+\infty$, while the trajectories of the other particles $(\lambda_1,\lambda_2,\dots,\lambda_{N-1})$ 
are extended in a continuous way (note again from \eqref{sde.ev} that for all $i\neq N$, $\lambda_i(t)$ has a limit when 
$t\uparrow \tau_1$). 
At each explosion $\tau_k$, we re-label the eigenvalues 
according to the circular change of indexation 
\begin{align}\label{new-index}
\lambda_1 \to \lambda_2 \to \dots \to \lambda_{N-1} \to \lambda_N \to \lambda_1\,. 
\end{align}

Now, in order to define the trajectory of the Hermitian process $(H(t))$ for all time $t\geq 0$, we need to 
check that the sequence
of explosion times $(\tau_k)_{k \geq 1}$, defined recursively for $k\geq 1$ as 
\begin{align*}
\tau_{k+1}:= \inf\{t \geq \tau_k: \lambda_N(t) = -\infty \}\,,
\end{align*}
 has no accumulation points in $\R_+$. This fact follows from 
\cite[Section 5]{bloemendal2} where the authors prove that the explosion times of the eigenvalues process \eqref{sde.ev}
correspond to the focal points, which are almost surely finitely many in compact sets of $\R_+$ (see in particular \cite[Proposition 5.1]{bloemendal2}).

In this way, the trajectory of the Hermitian process $(H(t))$ is defined  for all time $t\geq 0$.
Its eigenvalues process 
evolves according to the stochastic differential system \eqref{sde.ev} with a circular re labeling at each explosion time $\tau_k$, 
while the associated eigenvectors process follows \eqref{sde-vectors} with the same re labeling at time $\tau_k$.

Because $H(t)$ is invariant under rotation at all times, there is not much to say on the eigenvectors dynamics 
of the process $H$.  
In the next sections of this paper, we focus on the spectral statistics of $H$.

\begin{rem}\label{rem-norestart1}
One could have chosen a different dynamic where the process 
$(\lambda_1,\dots,\lambda_N)$ is still a solution of the stochastic differential system  
\eqref{sde.ev} but with no restarting procedure after the explosions. Instead the particles are killed in $-\infty$ whenever they explode at some finite time 
(in the sense that they are stuck forever in $-\infty$ and do not interact anymore with the living particles, which have not yet exploded). This interesting model seems more complicated to handle with our methods, see Remark \ref{rem-norestart} for more details.
\end{rem}

\section{Dynamics in the scaling limit} \label{dynamic-density}
We are mainly interested in the
empirical measure of the eigenvalues of the Hermitian process $H$ at time $t$
\begin{align}\label{emp-density}
\mu_t^N := \frac{1}{N} \sum_{i=1}^N \delta_{\lambda_i(t)}\,. 
\end{align}
Recall that the eigenvalues process $(\lambda_1(t),\lambda_2(t),\dots, \lambda_N(t))_{t\geq 0}$ satisfies 
the stochastic differential system  \eqref{sde.ev} on the intervals $(\tau_k;\tau_{k+1}), k \in \N$ with the restarting and re-indexing procedures \eqref{new-index} at each explosion time $\tau_k$.     

\subsection{Evolution equation for the spectral density}
We denote by $\mathcal{P}(\R)$ the space of probability measures on $\R$ and by 
$\mathcal{C}([0;T], \mathcal{P}(\R))$ the space  of continuous functions from 
$[0;T]\to \mathcal{P}(\R)$. 

Let us briefly recall the definition of the Stieltjes transform of a measure. 
Let $\H:= \{z\in \C: \Im z >0\}$ the upper half complex plane. 
If $\nu$ is 
a measure on $\R$, its Stieltjes transform\footnote{Note that the Stieltjes transform is sometimes defined as the negative of $G$ i.e. $\int_\R \nu(dx)/(z-x)$.} is the holomorphic function $G:\H\to \H$ defined 
by 
\begin{align*}
G(z) = \int_\R \frac{\nu(dx)}{x-z}\,. 
\end{align*}

We can recover the probability measure $\mu_t$ from 
the Stieltjes inversion formula, which writes for $x<y$, 
\begin{align}\label{inversion-formula2}
\lim_{\varepsilon \downarrow 0}\int_x^y \Im \, G(\lambda+i \, \varepsilon,t) \, d\lambda= \pi\, \mu_t[x;y] + \frac{\pi}{2}\big(\mu_t(\{y\}) - \mu_t(\{x\})\big)\,. 
\end{align}
Basic properties of Stieltjes transforms, which shall be useful throughout the paper, 
are recalled in Appendix \ref{stieltjes}. 

Our main result in this section establishes the convergence 
of the continuous stochastic process $(\mu_t^N)_{t\geq 0}$ 
when the dimension $N$ tends to $\infty$.

\begin{theorem}\label{main}
Let $T>0$ and $a \in \R$. 

Suppose that, at the initial time, the empirical density $\mu_0^N:= \frac{1}{N} \sum_{i=1}^N \delta_{\lambda_i(0)}$ converges 
weakly as $N$ goes to infinity towards some $\mu \in \mathcal{P}(\R)$.

Then,  $(\mu_t^N)_{0\leq t \leq T}$ converges almost surely in $\mathcal{C}([0,T],\mathcal{P}(\R))$\footnote{The space $\mathcal{C}([0,T],\mathcal{P}(\R))$ is a Polish space as $\mathcal{P}(\R)$ equipped with its weak topology is metrizable ($\R$ is a separable space)}. 
Its limit $(\mu_t)_{0\leq t\leq T}$ is the unique measure valued process such that 
$\mu_0=\mu$ and whose Stieltjes transform $G(z,t) :=\int_\R \frac{\mu_t(dx)}{x-z}$ 
satisfies the holomorphic equation 
\begin{align}\label{burgers-eq}
 G(z,t) = G(z,0) + \int_0^t  \partial_z  \left[   \frac{\beta}{4}  \, G(z,s)^2 + (z^2-a) \, G(z,s) + z  \right] \, ds \,. 
\end{align}
\end{theorem} 

The proof of Theorem \ref{main} is deferred to Section \ref{proof}.
Our approach is classical and follows the method introduced 
in \cite{agz,rogers} (see also \cite{jp-alice,satya}). 
It consists in writing an evolution equation for the Stieltjes transform $G_N(\cdot,t):\H\to \H$ 
of the probability measure $\mu_t^N$ thanks to It\^o's formula 
and the stochastic differential system \eqref{sde.ev} satisfied by the $\lambda_i(t),i=1,\dots,N$ 
in section \ref{evolution.GN}. We first prove the almost-sure pre-compactness of the family $((\mu_t^N))_{0\leq t \leq T},N\in \N$
in the space  $\mathcal{C}([0;T], \mathcal{M}_{\le 1}(\R))$ where $\mathcal{M}_{\le 1}(\R)$ 
is the space of measure with total mass smaller or equal to $1$ equipped with its weak-$\star$ topology. 
It turns out that it is sufficient for our purposes to prove pre-compactness in  $\mathcal{C}([0;T], \mathcal{M}_{\le 1}(\R))$ (instead of $\mathcal{C}([0;T], \mathcal{P}(\R))$)
as we can show afterwards that the Stieltjes transform solutions of \eqref{burgers-eq} 
are associated to probability measures. This feature is original and contrasts with the classical case of the quadratic potential for which the solutions of the Stieljes equation are not necessarily probability measures (see \cite{agz,rogers} where the authors prove a dynamical version of Wigner's Theorem). 
Finally we prove uniqueness of the solution \eqref{burgers-eq} in Lemma \ref{unicity-lemma}
using the characteristic method. 

\subsection{Convergence to equilibrium} 
We are interested in this paragraph in the convergence of the probability measure process $\mu_t$ when $t\to \infty$
(in the space of Radon measures endowed with the topology of weak convergence).  

To prove that $\mu_t$ converges weakly as $t\to \infty$ to some $\mu\in \mathcal{P}(\R)$, 
it is sufficient to show (see \cite{geronimo}) that the Stieltjes transform $G(\cdot,t)$ 
associated to $\mu_t$ converges point-wise for $z\in \H$ 
to the Stieltjes transform $G$ of the \emph{probability measure} $\mu$.  

Although we were not able to prove it, it is natural to expect that the following convergence holds:
\begin{conjecture}\label{conjecture1}
Let $(G(z,t))_{z\in \H,t\geq 0}$ be the solution of the evolution equation \eqref{burgers-eq}. 
Then, the following limit exists for all $z\in \H$, 
\begin{align}\label{def-G-infty}
G_a(z) := \lim_{t\to \infty} G(z,t) \,. 
\end{align}
\end{conjecture}
Assuming that conjecture \ref{conjecture1} is true, we can prove that the limit is indeed the Stieltjes transform of some probability measure and deduce the desired convergence:
\begin{proposition}\label{prop-conv-time}
The function $G_a$ is the Stieltjes transform of a probability measure $\mu_a$
and is a stationary solution of the evolution equation \eqref{burgers-eq}.

Consequently, the probability measure $\mu_t$ defined in Theorem \ref{main} converges weakly as $t\to +\infty$ 
to the probability measure $\mu_\infty$.  
\end{proposition}
The proof of this Proposition is deferred to Subsection \ref{proof-prop-conv-time}.

\medskip

In the next section, we characterize $G_a:\H \to \C$ as the unique function 
satisfying the following two properties:
\begin{itemize}
\item $G_a$ is a stationary solution of \eqref{burgers-eq};
\item $G_a$ is the Stieltjes transform of a (probability) measure $\mu_a$. 
\end{itemize}
Moreover, we derive $G_a$ and $\mu_a$ explicitly in analytic forms (see Theorem \ref{theorem-statio} and its proof below).
The probability measure $\mu_a$ is the limiting empirical eigenvalue density of the matrix $H(t)$ in the 
stationary state. We will see that $G_a$ provides additional information on the current of particles
in the system, when the well in $\sqrt{a}$ of the potential $V_a$ is not too confining to retain all the particles inside (see
section \ref{stationary-flux} for further details).

\section{Equilibrium spectral density}  \label{equilibrium-density}

In this section, we compute the limiting (i.e. when $N=\infty$) empirical eigenvalue density of the matrix $H(t)$ 
 in the stationary state (i.e. after a long time $t\to \infty$). 

Under the assumption that Conjecture \ref{conjecture1} is true, we know that there exists a Stieltjes
transform $G_a$ associated to $\mu_a\in \mathcal{P}(\R)$ which is a stationary solution of the evolution  
equation \eqref{burgers-eq}. 
In this section, we prove the uniqueness of the analytic function $G_a$ which enjoys those properties 
and we compute it explicitly in terms of the roots of a cubic polynomial. 
Using the Stieltjes inversion formula recalled in \eqref{inversion-formula}, 
we derive the probability density of $\mu_a$
with respect to Lebesgue measure and exhibit an interesting (sharp) phase transition displayed by 
 $\mu_a$
at the critical value $a=a^*=\frac{3}{4}\beta^{2/3}$.

Many details on the stationary probability measure $\mu_a$ are provided in this section. 
We  briefly summarize the main features in the following theorem. 

\begin{theorem}\label{theorem-statio}
For all $a \in \R$, there exists 
a unique analytic function $G_a$ with the following two properties: 
\begin{itemize}
\item it is the Stieltjes transform of a probability measure $\mu_a$, or equivalently  
(Akhiezer's Theorem \cite[page 93]{akhiezer}) $G_a$  is analytic on $\H$ with $G_a(\H)\subseteq \H$
and $ G_a(iy) \sim-1/(iy)$ as $y\to +\infty$; 
\item it is a stationary solution of the evolution equation \eqref{burgers-eq-statio}. 
\end{itemize}
  
Moreover, the probability measure $\mu_a(d\lambda)$ 
admits a density $\rho_a(\lambda)$ with respect to the Lebesgue measure, 
computed explicitly (see \eqref{subcritical-density} and \eqref{supercritical-density}) 
in terms of the roots of the polynomial of degree three $P'(z):= 4z^3-4a z -\beta$. 

We distinguish two regimes depending on whether $a\geq a^*$, or $a < a^*$ 
where 
\begin{equation*}
a^* = \frac{3}{4} \beta^{2/3}
\end{equation*}
 is the critical value at which the probability density $\rho_a$ displays a sharp phase transition 
 (illustrated in Fig. \ref{fig.potential}): 
\begin{itemize}
\item If $a\geq a^*$, $\rho_a$ is supported on a compact interval. 
\item If $a< a^*$, $\rho_a$ has full support in $\R$ and is flanked with 
symmetric heavy tails as $x\to \pm \infty$, 
\begin{align*}
\rho_a(x) \sim \frac{C_a}{x^2}\,,
\end{align*}
where  $C_a >0$ is an explicit constant (see below). 
\end{itemize} 
\end{theorem}

In the supercritical regime $a\geq a^*$, the particles are all confined in the well of the local minimum 
of the potential $V_a$. The well is deep enough compared to the electrostatic repulsion between the eigenvalues, 
to keep the particles inside it. 
The limiting density of particles has a classical shape 
in random matrix theory with a compact support and singularities of order $1/2$ at the 
upper and lower edges of the spectrum (square root cancelations). 
 The critical density $\rho_{a^*}$, still compactly supported, 
is of particular interest with an usual singularity at the lower edge (see paragraph \ref{critical-regime}
where we compute it explicitly as function of $\beta$ (see \eqref{density-critical}).

A sharp transition is observed at the critical value $a=a^*$: the 
 density $\rho_a$ is compactly supported for $a\geq a^*$
 but has full support with heavy tails if $a< a^*$. 
As we will see in the next section, the equilibrium 
 of the eigenvalues becomes very unusual in the subcritical 
 regime, mainly due to the non-confining shape of the potential $V_a$.
 This observation appears to be new.  
  Although the density profile of the eigenvalues has a stationary shape $\rho_a$,
 the particles are still flowing across the system in the equilibrium state. 
 There is in fact a positive current of eigenvalues
 flowing from $+\infty$ to $-\infty$ in a stationary way: the number of particles per unit of time 
 shifting from the right to the left (counted algebraically) 
 at some given level $x\in \R$ is constant (in time and in space). This 
 stationary current is also computed explicitly in the next section. 
 
\begin{figure}[h!btp] 
    \includegraphics[scale=0.85]{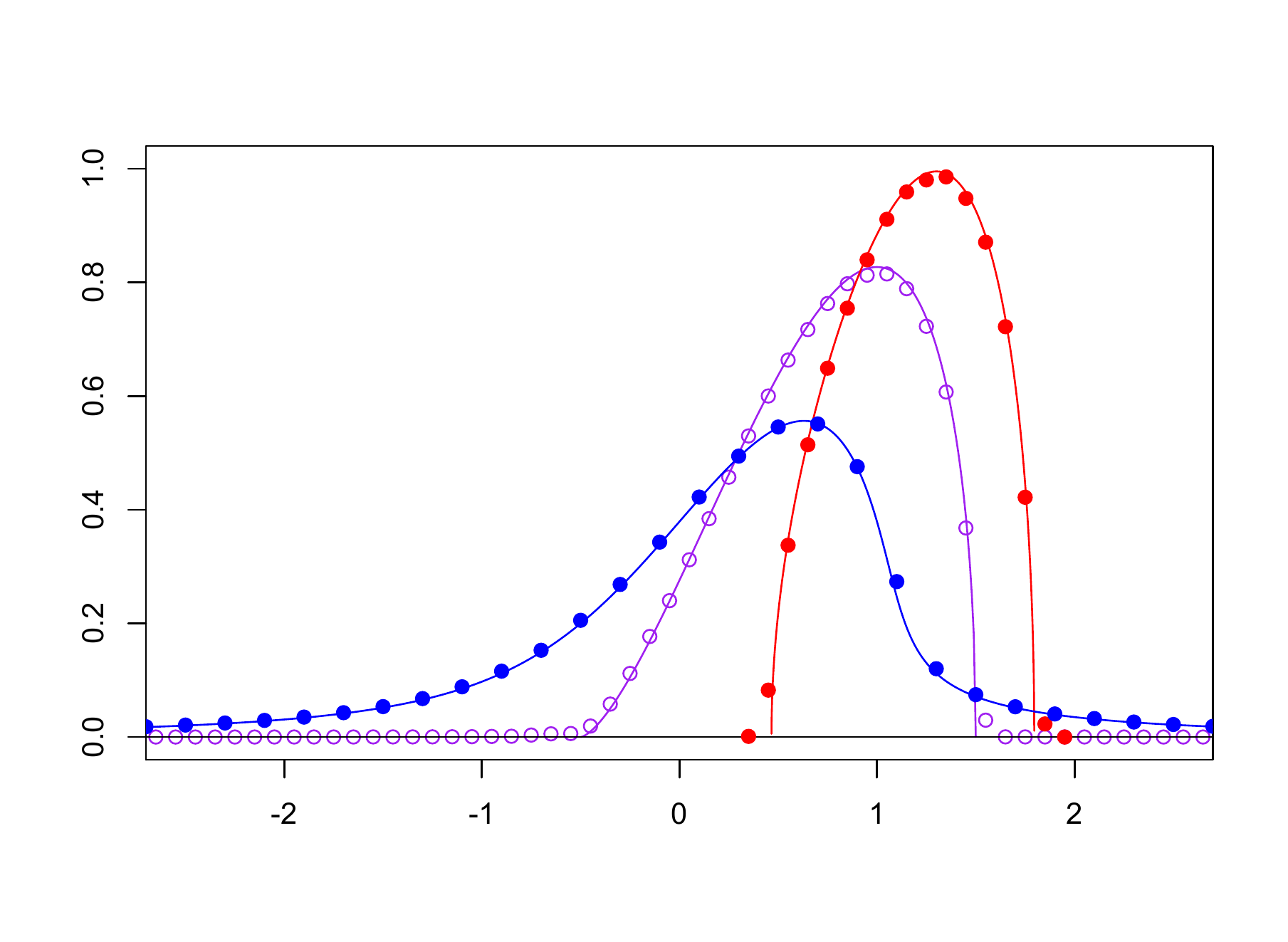}
     \caption{(Color online). Graphics 
     of the limiting eigenvalues densities $\rho_a(x)$ as a function of $x$, in the real case ($\beta=1$) for 
    $a=0$, $a=a^*= 3/4$ and $a= 3/2$ corresponding respectively to the sub-critical, 
    critical and super critical regimes. The points on top of the straight line represent the empirical densities obtained from our simulated samples with $N=50$. }\label{fig.eigenvalues-densities}
\end{figure}

We provide an illustration of the variety of possible behaviors for the density $\rho_a$ 
in Fig. \ref{fig.eigenvalues-densities}
where we show the graphics of the limiting eigenvalues density $\rho_a(x)$ as a function of $x$, 
for particular values of $a$ 
in the three different regimes $a<a^*, a= a^* $ and $a >a^*$. 
We have also checked our result with numerical simulations with excellent agreement 
(see Fig. \ref{fig.eigenvalues-densities}).  
 The samples to construct the empirical densities were obtained by simulating the
 Hermitian matrix process $(H(t))_{t\geq 0}$ satisfying \eqref{langevin-cubic}, with $N=50$. 
 The method is usual and consists in discretizing time and diagonalizing 
 the matrix $H(t)$ at each time step. We have introduced a cut off in order to deal with the explosions. 
 Whenever the lowest eigenvalue $\lambda_N(t)$ 
 of the matrix $H(t)$ gets smaller than the cut-off value, we re initialize
 the matrix $H(t)$ according to the procedure described in section \ref{cubic-hermitian-process}, using 
 again a cut-off to approximate the value $+\infty$. 
 We let our algorithm run for a time $t=100$ with a time 
 step $\delta t= 10^{-3}$ and
 constructed the empirical densities in the respective cases 
 using the eigenvalue samples at all time steps after time $10$. 
 We noted that both the convergences in dimension and time are extremely fast.

The rest of this section is devoted to the proof of Theorem \ref{theorem-statio}. 
Additional informations on the probability measure $\mu_a$ are provided along this proof.  
 
{ \it Proof of Theorem \ref{theorem-statio}}. 
From Conjecture \ref{conjecture1} and Proposition \ref{prop-conv-time}, we know that there exists 
an analytic function $G_a$ which satisfies the two properties given in Theorem \ref{theorem-statio}. 
It is therefore sufficient to prove the uniqueness of $G_a$ with those two properties. 

In the following, we work with the complex square root function $\sqrt{\cdot}: \C \to \H\cup \R^+$ 
defined for $z= r e^{i\theta} \in \C, r\geq 0, \theta \in [0;2\pi)$ as $\sqrt{z}:=\sqrt{r} \, e^{i\theta/2}$. 
With this definition, the square root is analytic on the domain $\C\setminus \R_+$ (to $\H$) but displays 
a discontinuity near the positive half line $\R_+$. 

If $G_a$ is a stationary solution of \eqref{burgers-eq}, then there exists a constant $J\in \C$ 
 such that for all $z\in \H$, 
 \begin{align}\label{burgers-eq-statio}
 \frac{\beta}{4}  \, G_a(z)^2 + (z^2-a) \, G_a(z) + z   = J \,. 
\end{align}

For any $z\in \H$, we can solve the second degree equation \eqref{burgers-eq-statio} 
for which we have two solutions, 
\begin{align}\label{solutions}
G_\pm(z) = \frac{2}{\beta} \left( a- z^2 \pm \sqrt{(z^2-a)^2 - \beta(z-J)} \right)\,. 
\end{align}
Our analytic function $G_a:\H\to \H$ is equal either  to $G_+(z)$ or to $G_-(z)$
depending on the value of $z\in \H$. It is possible that both $G_+(z)$ and $G_{-}(z)$ belong to $\H$ 
for certain values of $z$, but for any $z\in \H$, we will see that there will be  one unique possible value of $G_a(z)$.

We now seek for the 
constants $J$ such that $G_a$ has the two required properties of Theorem \ref{theorem-statio}. 
We shall in fact prove that 
there exists a unique such constant $J:=J_a$. 

The main idea is that the analyticity of $G_a:\H\to \H$ and the {\it non- analyticity} 
of the square root function in $0$ prevent the complex polynomial function $P(z):=  (z^2-a)^2 - \beta(z-J)$
inside the square root of \eqref{solutions} to have any root with {\it odd} multiplicity, one or three, in $\H$.

Additional information on the spatial locations of the 
roots of the polynomial $P$ is then provided by the zeroes of its derivative with degree three,
\begin{align}\label{Pprime}
P'(z) = 4 z^3 - 4 a z -\beta\,. 
\end{align}
Mainly, it is well known ({\it Gauss-Lucas} Theorem) that the roots of $P'$ 
all lie within the {\it convex hull} of the roots of $P$, that is the smallest convex polygon containing the roots of $P$.

It is easy to check from {\it Cardan}'s formulas that $P'$ has three real roots if and only if 
\begin{align*}
a \geq \frac{3}{4} \, \, \beta^{2/3}= a^*\,. 
\end{align*}
If $a< a^*$, then $P'$ has one real root and two roots in $\C\setminus \R$ which are complex conjugate.
Moreover Cardan's formulas permit us to compute the three roots of $P'$ analytically. Gathering those arguments, 
we can now compute the constant $J_a$ treating separately the two cases  $a< a^* $ and $a\geq a^* $.

\subsection{Subcritical regime}
If $a< a^* $, then $P'$ has one (unique) root in $\H$ so that, 
by Gauss-Lucas Theorem, $P$ has at least
one root in $\H$. By analyticity of $G_+$ on $\H$, any such root of $P$ is necessarily of multiplicity two \footnote{Multiplicity 
four is excluded because it would imply for example that $P'$ would have a root of multiplicity three.}, 
and therefore is equal to the unique root of $P'$ in $\H$. Denoting by $\zeta_a$ this root and using the condition $P(\zeta_a)=0$, we obtain 
\begin{align}\label{eq-J-alpha}
J_a &= \zeta_a -\frac{1}{\beta} (\zeta_a^2-a)^2\\ 
&\Big(=\frac{a}{\beta}\, \zeta_a^2 + \frac{3}{4}\, \zeta_a - \frac{a^2}{\beta}\Big)\,, \notag
\end{align}
and the polynomial $P(z)=P_a(z):= (z^2-a)^2 - \beta(z-J_a)$ is uniquely determined. 

Note that we have the following expression for $\zeta_a$:
\begin{align*}
\zeta_a = \frac{\beta^{1/3}}{2}\left(\Big(1+ \sqrt{1-(a/a^*)^{3}}\Big)^{1/3} 
\,\mathbf{j}\,+\Big{|}1- \sqrt{1-(a/a^*)^{3}}\Big{|}^{1/3}\, \mathbf{j^2}
\right)
\end{align*}
where $\mathbf{j} = -1/2+\sqrt{3} \,\mathbf{i}/2$ which permits us to check that $J_a \in \H$ thanks to elementary computations (using the second expression for $J_a$). 
Therefore $G_a$ can be extended analytically to a neighborhood of the real axis (indeed $\zeta_a \in \H$ implies that $P_a$ cannot have real roots). Moreover $G_a$ takes the following expression for $z\in \H$ near the real axis:
\begin{align}\label{final-formula-Ga-sub}
G_a(z) = \frac{2}{\beta}  \left( a- z^2 + \sqrt{P_a(z)}\right)\,. 
\end{align}

\begin{rem}
The polynomial $P_a$ can be factorized with elementary computations as 
\begin{align}\label{factorize-Pa}
P_a(z):= (z-\zeta_a)^2  (z-\gamma_-)(z-\gamma_+)\,,
\end{align}
where 
\begin{align*}
\gamma_\pm &= -\zeta_a \pm \sqrt{2(a-\zeta_a^2)} 
\in \H_{ -}:= \{z \in \C \;:\; \Im(z) < 0\}
\end{align*}
\end{rem}
 
Eq. \eqref{final-formula-Ga-sub} characterizes the function $G_a$ uniquely. 
The Stieltjes inversion formula \eqref{inversion-formula} permits us to check that 
the probability measure $\mu_a$ admits 
a density $\rho_a$ with respect to the Lebesgue measure given by for $x\in \R$,
\begin{align}\label{subcritical-density}
\rho_a(x) = \frac{2}{\beta\pi} \Im \left[\sqrt{P_a(x)} \right]\,. 
\end{align} 
It is easy to check that $\rho_a$ has full support in $\R$. Recalling that $\Im J_a >0$, it is straightforward 
to derive the {\it heavy} tails of $\rho_a$ when $x\to \pm \infty$, 
\begin{align}\label{tails}
\rho_a(x) \sim \frac{1}{\pi} \frac{\Im J_a}{x^2}\,. 
\end{align}  
The function $\rho_a$ is integrable on $\R$ as expected for a probability measure.      
It would be interesting to check the normalization condition 
$\int_\R \rho_a=1$ \footnote{We did a numerical check of this fact with mathematica.} by direct integration. 

In order to prove the existence of the Stieltjes transform $G_a$ 
satisfying the two required properties of Theorem \ref{theorem-statio}
{\it without} assuming that Conjecture \ref{conjecture1} holds, one would need to 
prove that the Stieltjes transform of the {\it explicit} probability density $\rho_a$ given in \eqref{subcritical-density}
is indeed the analytic function $G_a$ characterized in \eqref{final-formula-Ga-sub}.  
We were not able to perform this integration, although the two formulas \eqref{subcritical-density}
and \eqref{final-formula-Ga-sub} look very similar. 
In the super-critical regime (see below), we can do this integration.    

\subsection{Super critical regime}
If $a \geq a^*$, then the derivative polynomial $P'$ has three real roots, and 
the analyticity conditions on $G_a$ and the {\it Gauss-Lucas} Theorem
permit to show that all the roots of $P$ are real valued. 
Indeed, we have already seen that the polynomial $P$ can not have any root in 
$\H$ (otherwise, this root would be of multiplicity two and $P'$ would have a root in $\H$, which leads to a contradiction).
The remaining scenario where $P$ has two distinct zeroes in $\R$ 
with multiplicity one and two other zeroes (counting multiplicity) in $\H_-$ 
is also excluded: the zeroes of $P'$ would then 
lie on the frontier of the convex hull of the zeroes of $P$ and this is possible only if 
the two real zeroes of $P$ have multiplicity two, leading again to a contradiction.  

We conclude that all the zeroes of the polynomial $P$ are real, so that $P$ 
has real coefficients. In particular, $J \in \R$. The polynomial $P$ is now determined up 
to the real constant $J\in \R$, which has the effect of translating vertically 
(along the $y$-axis) the graph $\{(x,P(x)):x\in \R\}$.

The uniqueness of $J:=J_a\in \R$ leading to the correct solution $G_a$
will come from the normalization condition $\int_\R \mu_a(dx) =1$.    
By the Stieltjes inversion formula \eqref{inversion-formula}, 
we know that the measure $\mu_a$ is supported on the compact set
$K:=\{ x\in \R: P(x) \leq 0\}$ \footnote{$K$ is a union of intervals. For instance, if $P$ has four distinct eigenvalues $x_1<x_2<x_3<x_4$, 
then $\mu$ has a disconnected support of the form $[x_1;x_2] \cup [x_3;x_4]$.} 
and has a density $\rho_a$ with respect to Lebesgue, defined for $x\in K$ as,
\begin{align*}
\rho_a(x) = \frac{2}{\beta\pi} \, \sqrt{-P(x)} =  \frac{2}{\beta\pi} \, \sqrt{\beta(x-J) -(x^2-a)^2 } \,. 
\end{align*}  
The area $\mathcal{A}$ 
under the graph of $\rho_a$ can thus be seen as a function of $J$, $\mathcal{A}:=\mathcal{A}(J)$. 
It is obviously a continuous and strictly decreasing function of $J$, and we have 
$\lim_{J\to -\infty} \mathcal{A}(J) = +\infty$
and $\lim_{J\to +\infty} \mathcal{A}(J) = 0$ (for $J$ large enough, $P(x)>0$ for all $x\in \R$). 
By the intermediate value theorem, there exists a unique $J:=J_a$ such that $\mathcal{A}(J)=1$. The uniqueness
of $J:=J_a$ is proved.  

Now we would like to determine the value of $J_a$ in the present case, $a \geq  a^*$. To guess its value, let us notice that
if $P$ has four distinct real roots in $\R$, 
then the measure $\mu_a$ 
has a disconnected compact support, union of two disjoint intervals, and this solution is not physically sound
in view of the shape of the potential $V_a$. 
Recalling the confining shape of the potential $V_a(x)$ 
near the region $x\sim \sqrt{a}$ when $a > 0$, we would rather expect the polynomial $P$ 
to have its smallest root of multiplicity two and then two other roots of multiplicity one
near the confining zone of the potential. 
The minimal root of $P$, which will be denoted (again) by $\zeta_a$, would then also be the 
minimal root of $P'$.  
Finally, we can compute the real constant $J:=J_a$ associated to this scenario 
using the condition $P(\zeta_a)=0$, 
and we re obtain, with now $\zeta_a\in\R$,  formulas \eqref{eq-J-alpha} for $J_a$, \eqref{factorize-Pa}  
for the polynomial $P_a$
and \eqref{final-formula-Ga-sub} for $G_a(z)$ when $z$ is near the real line.

From the Stieltjes inversion formula, we see that $\mu_a$ has a density 
$\rho_a$ with compact support $[\gamma_-;\gamma_+]$, given by 
\begin{align}\label{supercritical-density}
\rho_a(x) 
=  \frac{2}{\beta\pi} (x-\zeta_a) \sqrt{(x-\gamma_-)(\gamma_+-x)}
\end{align}  
for $x\in [\gamma_-;\gamma_+]$. 
Reciprocally, we can check with an elementary integration that $\rho_a$ is 
 a probability measure and (using the residue Theorem) 
that its Stieltjes transform is indeed the analytic function $G_a$ 
characterized for $z$ near the real line in \eqref{final-formula-Ga-sub}.
We therefore do not need to assume Conjecture \ref{conjecture1} to prove the existence of $G_a$ in this case.  We treat in the next paragraph the particular case when $a=a^*$ in more details.

\subsection{Critical regime }\label{critical-regime}

We now analyze (even more explicitly) the probability density of the measure $\mu_{a^*}$ at the critical value. 
A straightforward computation permits to factorize the derivative polynomial as
\begin{equation*}
P'(z) = 4 \, (z+\frac{\beta^{1/3}}{2})^{2} \, (z-\beta^{1/3})\,. 
\end{equation*}
Its minimal root (with multiplicity two) is $\zeta_{a^*}= - \frac{1}{2} \, \beta^{1/3}$ and from relation \eqref{eq-J-alpha}, we find \begin{equation*}J_{a^*}=-\frac{3}{4}\beta^{1/3}.\end{equation*} 

The factorization of $P_a$ writes as    
\begin{align*}
P_a(z) = (z+ \frac{\beta^{1/3}}{2})^3 \, (z-\frac{3}{2} \beta^{1/3})\,. 
\end{align*}
We note in particular that $\zeta_a$ is of multiplicity three for $P$. This behavior is rather natural at the transition:
the root $\zeta_a$ with multiplicity two for $a>a^*$ reaches multiplicity three for $a=a^*$ (the first and 
second roots merge together) and finally splits up into three non real roots in $\C\setminus \R$, leading to a constant $J_a\in \H$. 
We shall see in the next subsection that the probability density starts to flow at this point. 

We finally derive the density $\rho_{a^*}$, which has a compact support,
\begin{align}\label{density-critical}
\rho_{a^*}(x) = \frac{2}{\beta\pi} \, \left(x+  \frac{1}{2} \, \beta^{1/3}\right)^{3/2} \, \sqrt{\frac{3}{2} \beta^{1/3}-x}\,, 
\quad   - \frac{1}{2} \, \beta^{1/3} \leq  x \leq \frac{3}{2} \beta^{1/3}\,. 
\end{align}
We note the unusual order of the singularity $3/2$ 
of the density $\rho_{a^*}$ at its lower edge $- \frac{1}{2} \, \beta^{1/3}$.  
At the critical value, the concavity changes near the lower edge in order to prepare the 
sharp phase transition when $a$ gets smaller than $a^*$.
 The probability will start to flow in the system as soon as $a < a^*$
with a constant profile (allocation) $\mu_a$  supported on the whole real line $\R$ flanked with heavy tails.   
It is interesting to note that such singularities at the edges of random matrices spectrum 
were already observed in \cite{bipz} in the context of a non-confining quartic potential. 
We shall revisit some of the questions investigated in \cite{bipz} in section \ref{section-quartic-potential}. 

We are now interested in the so called flux of particles in the system, which 
measures the number of particles which shifts from the right to the left at some given level $x\in \R$ 
per unit of time. It turns out that our method permits us to compute this flux explicitly in the stationary regime. 

\section{Stationary flux of charges}\label{stationary-flux}
Let us recall the evolution equation \eqref{burgers-eq} which may be rewritten for $s< t$ as 
\begin{align}\label{eq-G-J}
G(z,t) - G(z,s) = \int_s^t \partial_z \, J(z,u) \, du\,,
\end{align}
where 
\begin{align}\label{def-J}
J(z,u) = \frac{\beta}{4}  \, G(z,u)^2 + (z^2-a) \, G(z,u) + z\,. 
\end{align}
Eq. \eqref{eq-G-J} is a continuity equation. We can transform \eqref{eq-G-J} into an evolution equation 
on the measure $\mu_t$ thanks to the Stieltjes inversion formula recalled in \eqref{inversion-formula2}, so that the 
interpretation of this continuity equation becomes clearer. 
Taking imaginary 
part, integrating \eqref{eq-G-J} on the horizontal segment $[x+i\varepsilon;y+i\varepsilon], x<y, \varepsilon >0$ and sending 
$\varepsilon \to 0$, we obtain 
\begin{align*}
 \mu_t[x;y] - \mu_s[x;y] =\frac{1}{\pi} \lim_{\varepsilon\to 0} \int_{x+i\varepsilon}^{y+i\varepsilon} dz \int_s^t \partial_z \, \Im J(z,u) \, du\,. 
\end{align*}
For $\varepsilon >0$ fixed, the Fubini Theorem permits us to exchange the order of integration over $z$ and $u$ in the right hand side. 
We eventually obtain, for any $s<t, x<y$, 
\begin{align}\label{continuity-eq-mu}
 \mu_t[x;y] - \mu_s[x;y] =\frac{1}{\pi} \lim_{\varepsilon\to 0} \int_s^t \Im \,J(y+i\varepsilon,u) \, du  -  \int_s^t \Im \, J(x+i\varepsilon,u)  \, du \,. 
\end{align} 
We can interpret the probability measure $\mu_t$ as the electrostatic charge flowing across $\R$ from $+\infty$ to $-\infty$. Therefore, the right hand side of \eqref{continuity-eq-mu} may be seen as the amount of charge which enters the interval $[x;y]$ during the time interval $[s,t]$.

In order to further extend the present discussion, we 
admit that the Stieltjes transform $G(z,t)$ solution of \eqref{burgers-eq} has a \emph{continuous extension to $\R \cup \H$}. Unfortunately, we are not able to prove this mathematical detail, although it is physically sound. 
Note that this continuous extension was proved in \cite{biane} in the case of the complex Burgers equation given in 
\cite[Introduction]{biane}, which is the free analogue of the heat equation. 
The fact that $G(z,t)$ has a continuous extension to $\R \cup \H$ implies that the probability measure $\mu_t(dx)$ admits 
a density $\rho_t(x)$ with respect to Lebesgue measure, such that $(x,t) \mapsto \rho_t(x)$ is smooth, 
and we have 
\begin{align*}
\lim_{\varepsilon \downarrow 0} G(\lambda+ i \varepsilon,t) = P.V. \int_\R \frac{\rho_t(x) }{x-\lambda} dx + i \pi \rho_t(\lambda)\,,
\end{align*}
 where $P.V.$ stands for Principal Value. 

Under this assumption, it is clear that the analytic function $J(\cdot,t)$ defined in \eqref{def-J} 
has a continuous extension to $\R \cup \H$. 
In particular, we have
\begin{align}\label{conv-Im-J-real-axis}
\lim_{\varepsilon \downarrow 0} \Im \, J(\lambda+ i \varepsilon,t) = \pi j_t(\lambda)\,, 
\end{align}
where 
\begin{align*}
j_t(\lambda) =  \rho_t(\lambda) \left(\frac{\beta}{2} P.V. \int_\R \frac{\rho_t(x)}{x-\lambda} dx + \lambda^2 -a \right)\,. 
\end{align*}
Coming back to \eqref{continuity-eq-mu}, we easily obtain using \eqref{conv-Im-J-real-axis}
\begin{align}\label{continuity-eq-mu-j}
 \mu_t[x,y] - \mu_s[x,y] =  \int_s^t (j_u(y) - j_u(x))\, du\,,
\end{align}
which may be rewritten also as $\partial_t \int_x^y \rho_t(\lambda) d\lambda = \frac{1}{\pi} (j_t(y)- j_t(x))$. 
From \eqref{continuity-eq-mu-j}, we have a clear physical interpretation of the quantity $j_t(x)$ which is precisely 
the flux of probability density in $x$ at time $t$, measuring the amount of probability density shifting 
from the right to the left of $x$ (algebraically) per unit of time at time $t$. 

The interesting feature of our matrix model \eqref{langevin-cubic} is that 
the flux does not vanish identically in the stationary state of the scaling limit, as we will see.

If Conjecture \ref{conjecture1} holds such that $G(z,t)\to_{t\to \infty} G_a(z)$ pointwise in $\H$, 
then it follows immediately that 
$J(z,t) \to_{t\to +\infty} J_a$.
Using now \eqref{conv-Im-J-real-axis} and Montel's theorem, we get 
\begin{align}\label{conv-j-infty}
j_t(x) \longrightarrow  _{t\to\infty} \frac{1}{\pi} \, \Im J_a\,. 
\end{align}
Note that, as one may have expected, 
the flux of probability becomes independent of the position $x\in \R$ in the stationary state. 
Recall also from section \ref{equilibrium-density} that $\Im J_a$ is non 
zero if and only if $a<a^*$. In fact, if $a\geq a^*$, the well of the potential $V_a$ is deep enough to confine 
all the particles. The number of explosions per unit of time in the stationary state is negligible
compared to the macroscopic mass of the confined particles. On the other hand, if $a < a^*$, then the well 
in $\lambda=\sqrt{a}$ (see Fig. \ref{fig.potential}) is not strong enough to confine all the charges 
which repel each other with electrostatic interaction.    

If $a<a^*$, we have the analytic expression \eqref{eq-J-alpha} 
for $J_a$ in terms of the unique root $\zeta_a\in\H $ of $P'$, 
given in \eqref{Pprime}.  

We end this section by establishing the link with the discrete setting where the $N$ eigenvalues of the matrix process $(H(t))$
defined in \eqref{langevin-cubic}
are diffusing in the potential $V_a$ (see fig. \ref{fig.potential}).
In this context, the probability (or electrostatic charge for a physical analogy) is carried by the eigenvalues $\lambda_i(t)$ satisfying 
the stochastic differential system
\eqref{sde.ev}. Each particle $\lambda_i$ carries a proportion $1/N$ of the total probability. 
Using \eqref{conv-j-infty}, we conclude that in the large $N$ limit and in the stationary state (i.e. after a long time $t$), 
the flux of probability density is $\Im J_a/\pi$. In other words, the numbers of particles $N_t(x)$
which shift in $x$ from the right to the left at time $t$ per unit of time is, in the stationary state, proportional to $N$ with 
\begin{align}\label{leading-flux}
N_t(x) \sim \frac{N}{\pi} \, \Im J_a  \,. 
\end{align}
This formula was checked numerically with very good agreement for $N=50$. 
We have also noted that the convergence in time and in $N$ is very fast. 
It would be very interesting to compute the fluctuations of this flux of particles around its typical value. 
We leave this challenging problem for future research. 

\begin{rem}\label{rem-norestart}
Let us say a few words about the non-conservative system already mentioned in Remark \ref{rem-norestart1}, where the particles are killed when they explodes instead of being restarted at $+\infty$. One can easily adapt the proofs of Theorem \ref{main}, and prove that for any $T>0$ and $a \in \R$, if the empirical density at the initial time converges weakly as $N$ goes to infinity towards some $\mu \in \mathcal{P}(\R)$, then the limit points $(\mu_t)_{0\leq t\leq T}$ of the (almost surely pre-compact) sequence $(\mu_t^N)_{0\leq t \leq T}$ $\in \mathcal{C}([0,T],\mathcal{M}_{\le 1}(\R))$
satisfy $\mu_0=\mu$ and their Stieltjes transforms $G(z,t)$ 
are solution of the holomorphic equation 
\begin{align}\label{burgers-non-conservative}
 G(z,t) = G(z,0) + \int_0^t  \partial_z  \left[   \frac{\beta}{4}  \, G(z,s)^2 + (z^2-a) \, G(z,s) + \mu_s(\R) z  \right] \, ds \,. 
\end{align}

 The main difference with Theorem \ref{main} is that the solutions $(\mu_t)$ of \eqref{burgers-non-conservative}
 evolve in the space of measures with a total mass $\mu_t(\R)$ decreasing over time. 
 There is no longer uniqueness of a solution $(\mu_t)_{0\leq t\leq T}$ to  \eqref{burgers-non-conservative} such that $\mu_0=\mu\in \mathcal{P}(\R)$ (Eq. \eqref{burgers-non-conservative} depends itself of $\mu_{\cdot}(\R)$) and Eq. \eqref{burgers-non-conservative} does not characterize the limiting process $(\mu_t)_{0\leq t \leq T}$. It makes the analysis of this model much more complex than in the restarting case. 

We conjecture that the limit $t \to \infty$ of the continuum process $(\mu_t)$ corresponds to the \emph{metastable} equilibrium of the finite-$N$ model (note that for any finite $N$, all the particles will explode in a finite time almost surely). In the case of the sub-critical regime ($a < a^*$), the system should lose some mass until it gets closer and closer to the critical regime. The study of the stationary regime gives the relation $a_c := (3/4)\, \alpha^{1/3} \beta^{1/3}$ between the critical point $a_c$ and the mass $\alpha$ of the stationary measure.
In the supercritical regime, the behaviour of the system should not depend on the choice of restarting or killing the particles as the explosions are too scarce to matter.
\end{rem}

\section{Non-confining Quartic potential}\label{section-quartic-potential}
We now revisit a problem investigated in \cite{bipz} (see also \cite{douglas,biane2}) 
and related to the quartic potential $U_g$ such that
\begin{align*}
U_g(x) = \frac{x^2}{2} + g\, x^4\,, 
\end{align*}
where $g\in \R$ is usually referred as the coupling constant. For negative values of $g$, the potential $U_g$ is non-confining with 
$U_g(x)\to \pm \infty$ as $x\to \pm \infty$, respectively. In \cite{bipz}, the authors consider
an ensemble of random matrices with invariant law in the space of complex Hermitian matrices given by
\begin{align}\label{ensembles-bipz}
P(dH) = \frac{1}{Z} \exp\left(- \frac{1}{2} {\mbox{Tr}} H^2 - \frac{g}{N} \, {\mbox{Tr}} H^4 \right) \ dH\,.
\end{align}
Although the probability distribution \eqref{ensembles-bipz} does not make sense for $g<0$,  
the authors of \cite{bipz} are still able to derive (analytically) a density, which corresponds if $g>0$ to the 
limiting spectral density when $N\to \infty$ of the random matrices in the ensemble \eqref{ensembles-bipz}.  
 The probability density they obtain still makes sense even for $g\in [-\frac{1}{48}; 0)$. 
 The purpose of this section is to bring new lights on this result. We extend (with more details) 
 their computations for $\beta=2$ to general values of $\beta>0$.

\begin{figure}[h!btp] 
    \includegraphics[scale=0.85]{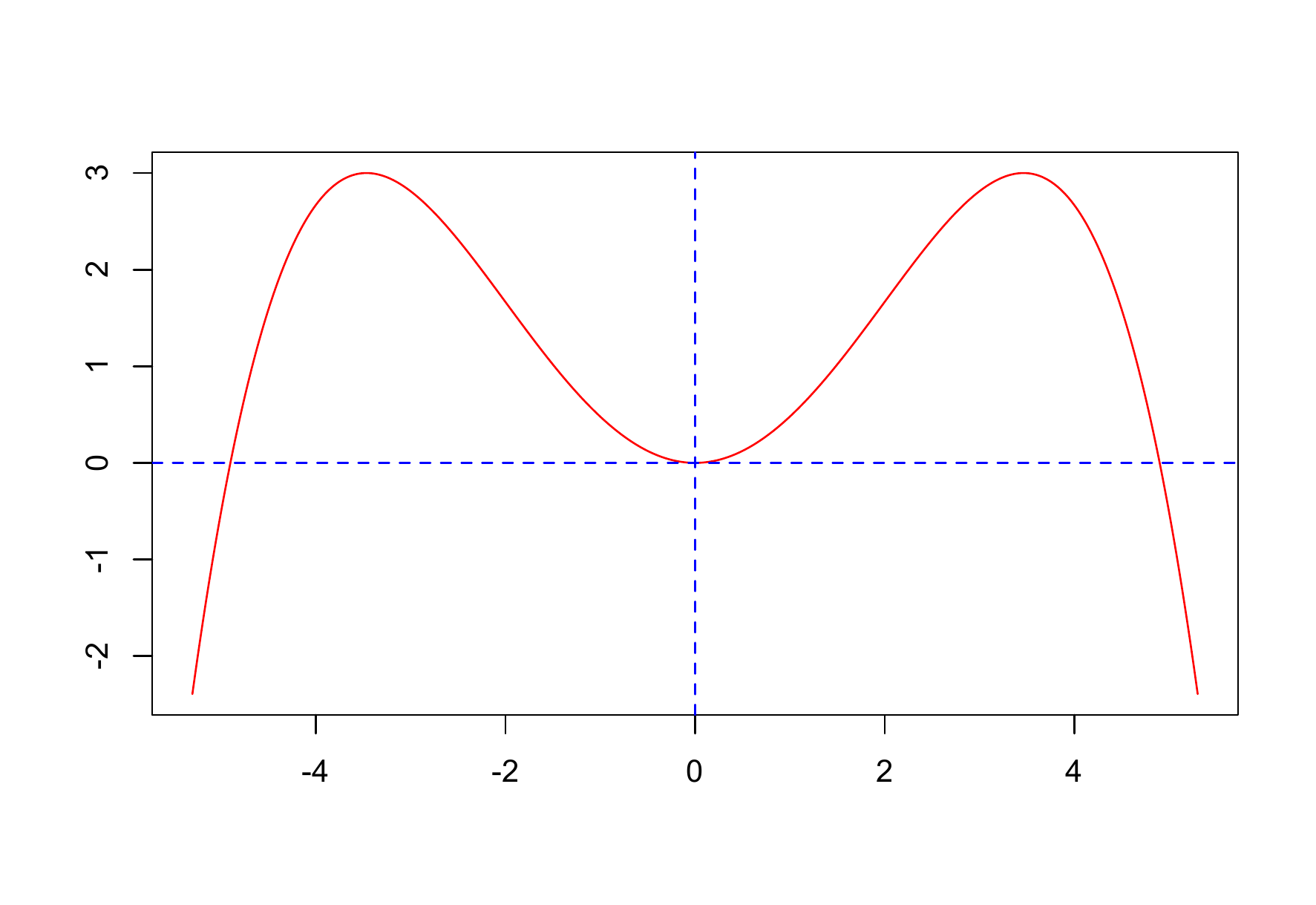}
     \caption{The potential $U_g(x)$ as a function of $x\in \R$ for $g=-1/48$. }\label{fig.potential-bipz}
\end{figure}

\begin{figure}[h!btp] 
    \includegraphics[scale=0.9]{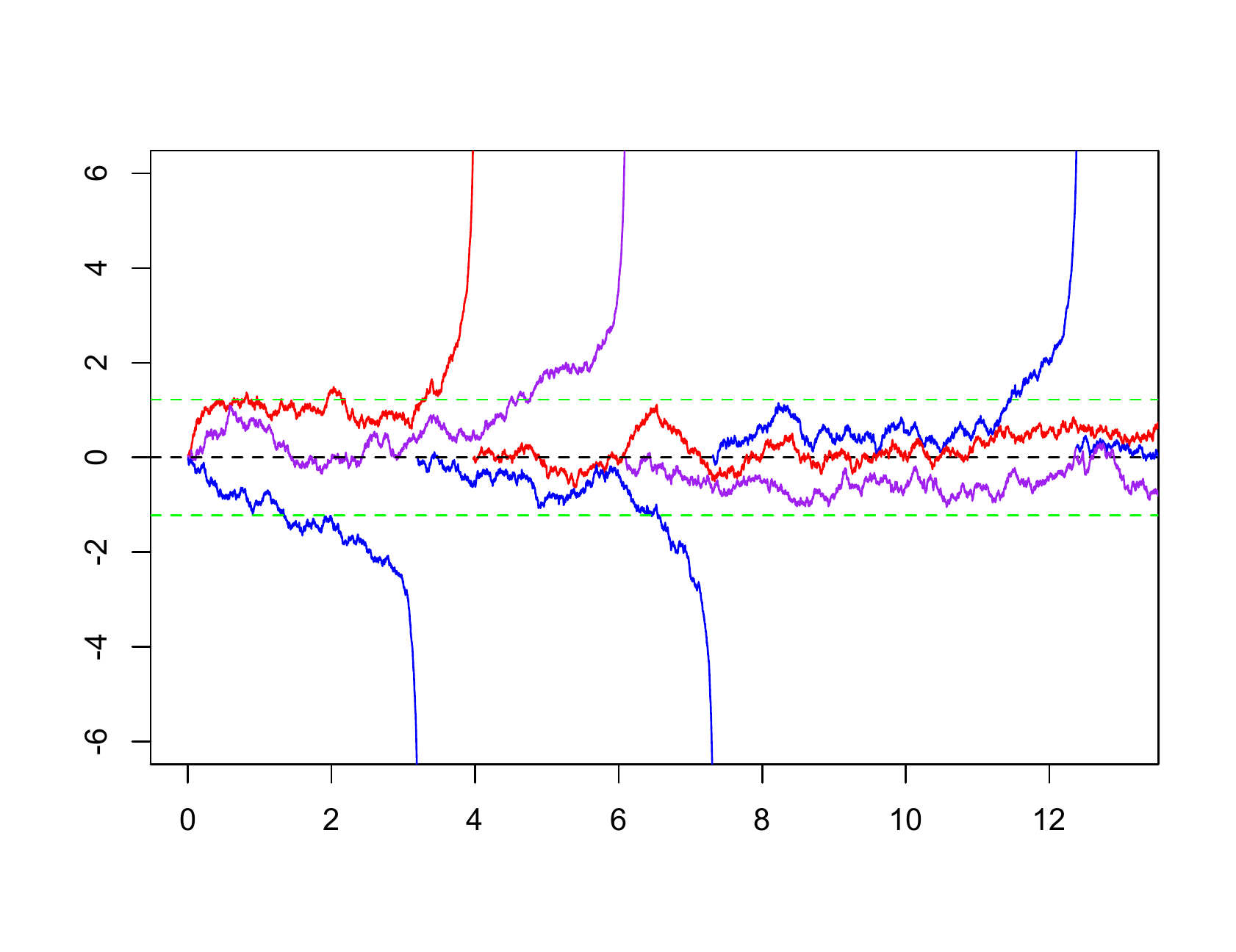}
     \caption{(Color online). Simulated paths of the eigenvalues $(\lambda_1(t),\lambda_2(t),\lambda_3(t))$ 
      as a function of time $t$ 
     for $N=3$, $\beta=1$ and $g=-1/6$. 
     The green horizontal dashed lines give the position of the  two symmetric local maximum (hills) of the potential 
     $U_g$ in $\lambda=\pm \sqrt{-1/(4 g)}$. }\label{fig.eigenvalues-path-Quartic}
\end{figure}

The method described above with the non-confining cubic potential $V_a$ permits 
us to define a stationary ensemble of random matrices in the potential $U_g$ for negative values of $g$. 
As before, we restrict ourselves to symmetric real random matrices although the complex and quaternion Hermitian cases
may be covered as well.
The idea is again to consider the symmetric matrix process $(H(t))$ such that 
\begin{align}\label{langevin-quartic}
dH(t) = - \left(\frac{1}{2} \, H(t)+ 2 g \, H(t)^3\right)\,  dt + \frac{1}{\sqrt{N}}\, dB(t)\,,
\end{align}
with $H(0)=0$ at the initial time and 
where $(B(t))$ is a  $N\times N$ symmetric Brownian motion, until the first explosion time.  
To extend the trajectory of the Hermitian process $(H(t))$ after this explosion time, we follow the method explained 
in section \ref{cubic-hermitian-process}. The situation is very similar: the non exploding eigenvalues and eigenvectors trajectories are a.s. continuous at the 
explosion times. The exploding eigenvalue is immediately restarted from $0$ (instead 
of $+\infty$ in section \ref{cubic-hermitian-process}). 
There are now two different types of explosions, either on the right side in $+\infty$ or on the left side in $-\infty$. 
Re starting the eigenvalues in $0$ permits us to preserve the symmetry with respect to $x=0$ and,  
for any $t\geq 0$, the equality in law $H(t)\stackrel{(d)}{=} - H(t)$ holds. 
A sample path of the eigenvalues of the process $H$ is shown in Fig. \ref{fig.eigenvalues-path-Quartic}. 
It is obtained with numerical simulations of the process $H$ following the restarting procedure at each explosion times described in
this paragraph. 

The restarting procedure of the eigenvalues in $0$ leads to some new difficulties, 
compared to the previous case. 
Mainly, the Stieltjes transform $G_N(z,t)$ of the empirical measure of the 
eigenvalues is not continuous at the explosion times. It is in fact right continuous in $\tau_k+$ with left limit in $\tau_k-$,  
with a jump of size $-\frac{1}{Nz}$ in $\tau_k$. 
In order to obtain an evolution equation for the empirical measure in a differential form, 
we can consider the set of smooth test functions  
\begin{align*}
 \mathcal{F}:= \{ f\in\mathcal{C}_0(\R): f(0)=0, \quad  x^3 \, f'(x)   \mbox{  is bounded on } \R\} \,. 
\end{align*}
 If $f\in \mathcal{F}$, then the function $t\mapsto \int_\R f(\lambda) \mu_t^N(d\lambda)$ is continuous on $\R_+$, even at the explosions times,
and we can use It\^o's formula to obtain, for any $t \geq 0$,  
\begin{align*}
d  \int_\R f(\lambda) \mu_t^N(d\lambda) &= - \int_\R \left(\frac{\lambda}{2}+2 g \lambda^3 \right) f'(\lambda) \mu_t^N(d\lambda) dt
+ \frac{\beta}{4} \int_\R \int_\R \frac{f'(\lambda)-f'(\lambda')}{\lambda-\lambda'} \mu_t^N(d\lambda) \mu_t^N (d\lambda') dt \\
&+\frac{1}{2N}(1-\frac{\beta}{2}) \int_\R f''(\lambda) \mu_t^N(d\lambda) dt +\frac{1}{N^{3/2}} \sum_{i=1}^N f'(\lambda_i) dB_i(t)\,,
\end{align*} 
where the $B_i$ are independent real Brownian motions. In order to prove convergence as $N \to +\infty$ of the empirical measure 
$\mu_t^N$ towards a probability measure $\mu_t$ and to characterize the evolution of the process $(\mu_t)_{t\geq 0}$, 
we propose to follow the same steps as in Section \ref{proof}. First, we show precompactness of the family 
of continuous processes
$(\mu_t^N)_{0\leq t\leq T}, N\in \N$ in the space $\mathcal{C}([0,T],\mathcal{M}_{\le 1}(\R))$. This 
proof is rather straightforward adapting the proof of Lemma \ref{precompactness-lemma} to the present quartic case. 
Then, taking $N\to \infty$, we obtain, as in subsection \ref{evolution.GN}, the following evolution equation for 
the limiting process $(\mu_t)_{0\leq t \leq T}$: 
for all $t\in [0,T]$ and $f\in \mathcal{F}$,
\begin{align}\label{evolution-mu-quartic}
\partial_t \int_\R f(\lambda) \mu_t(d\lambda) = - \int_\R \left(\frac{\lambda}{2}+2 g \lambda^3 \right) f'(\lambda) \mu_t(d\lambda) 
+ \frac{\beta}{4} \int_\R \int_\R \frac{f'(\lambda)-f'(\lambda')}{\lambda-\lambda'} \mu_t(d\lambda) \mu_t (d\lambda') \,. 
\end{align}
It would remain to show the uniqueness of the continuous process $(\mu_t)_{0\leq t\leq T}$ and that this process takes values
in the space of probability measures.  As in the proof of Theorem \ref{main}, we only know \emph{a priori} that the limit point 
(along some subsequence) 
$(\mu_t)_{0\leq t\leq T}$ belongs to  
$\mathcal{C}([0,T],\mathcal{M}_{\le 1}(\R))$ (instead of $\mathcal{C}([0,T],\mathcal{P}(\R))$). This question seems more difficult and we do not have insights on the method 
to prove this. 
We leave this problem as an open question. 

We are now interested in the probability measures $\mu$ which are stationary solutions 
of the evolution equation \eqref{evolution-mu-quartic}, i.e. such that 
for any test functions $f\in \mathcal{F}$,
\begin{align}\label{equilibrium-mu-quartic}
\int_\R \left(\frac{\lambda}{2}+2 g \lambda^3 \right) f'(\lambda) \mu(d\lambda) 
= \frac{\beta}{4} \int_\R \int_\R \frac{f'(\lambda)-f'(\lambda')}{\lambda-\lambda'} \mu(d\lambda) \mu(d\lambda')\,. 
\end{align}
Solving equation \eqref{equilibrium-mu-quartic} by choosing a particular family of functions $f\in \mathcal{F}$ 
 is not an easy task and has led us to heavy computations.  
Another route like in \cite{bipz} is to fix $z\in \H$ and simply consider $f(\lambda)=1/(\lambda-z)$ 
 (although $f\not\in\mathcal{F}$).
Equation \eqref{equilibrium-mu-quartic} conveniently rewrites (after a further integration) in term of the Stieltjes transform $G$ of the probability measure $\mu$ as 
\begin{align}\label{Stieltjes-quartic-equilibrium}
\frac{\beta}{4} G^2(z) + \left(2 g z^3 + \frac{z}{2}\right) G(z) + 2 g z^2 = J\,,
\end{align} 
where $J\in \C$ is an integration constant. 

We repeat the same steps as in section \eqref{equilibrium-density} 
to prove the uniqueness of the analytic function $G_g:\H\to \H$ with the following two properties: 
\begin{itemize}
\item there exists $J \in \C$ such that $G_g$ satisfies the 
quadratic equation \eqref{Stieltjes-quartic-equilibrium} for all $z\in \H$; 
\item $G_g$ is the Stieltjes transform of a probability measure $\mu_g$.
\end{itemize} 
In fact, assuming the existence of such a function $G_g$, we show that there is a unique 
possible constant $J:=J_g$ that we compute explicitly. We can then determine $G_g$ uniquely by solving the 
quadratic equation \eqref{Stieltjes-quartic-equilibrium}. Note however that, in order 
to prove existence,  
we still have  to check that 
the solution $G_g$ we obtain is indeed the Stieltjes transform of a probability measure. 
As we will see, such a Stieltjes transform $G_g$ exists if and only if $g\geq g_c  :=-\frac{1}{24\beta}$.  
For $g< g_c$ however, the method breaks down, leading to a 
measure which is not a probability measure. 
In this case, there are no Stieltjes transforms
satisfying \eqref{Stieltjes-quartic-equilibrium} whatever the values of $J\in \C$. 
For $\beta=2$, we re obtain the critical value $-1/48$ already obtained in \cite{bipz}.

We now explain the main steps of this computation. 
For any $J\in \C$ and all $z\in \H$, there are two solutions to the quadratic equation \eqref{Stieltjes-quartic-equilibrium}
\begin{align*}
G_\pm(z) = \frac{2}{\beta} \left(-\frac{z}{2} - 2 g z^3 \pm \sqrt{(2 g z^3 +\frac{z}{2})^2- \beta (2 g z^2 -J)} \right)\,. 
\end{align*}

The idea is again to compute the derivative 
polynomial $P_g'$ of the discriminant of the quadratic equation \eqref{Stieltjes-quartic-equilibrium}, which writes as 
\begin{align*}
P_g'(z) = z \left(24 g^2 z^4 + 8 g z^2 + \frac{1}{2} - 4 \beta g\right) \,.  
\end{align*} 
The quadratic equation $24 g^2 X^2 + 8 g X + \frac{1}{2} - 4 \beta g =0$ has two roots 
\begin{align*}
X_{\pm} = -\frac{1}{6 g} \left(1\pm\frac{1}{2} \sqrt{1+ 24 \beta g}\right)\,.
\end{align*}
Therefore, if $g> g_c$, all the zeroes of the polynomial $P'_g$ are real and we denote by $\xi_g$ the minimal root of $P_g'$, 
given for $g<0$, by
\begin{align*}
\xi_g:= - \sqrt{X_+} = \left(-\frac{1}{6 g} \left(1+ \frac{1}{2}\sqrt{1+ 24 \beta g}\right)\right)^{1/2}\,.  
\end{align*}
Using the Gauss Lucas Theorem and the analyticity of $G_g$, we can prove as before 
that all the zeroes of the polynomial $P_g(z):=  (2 g z^3 +\frac{z}{2})^2- \beta (2 g z^2 -J)$ are also real.  
This implies that $J\in \R$. With the same continuity argument of the area below the graph 
of the underlying probability density with respect to $J$, we can prove that there is in fact 
a unique $J_g\in \R$ insuring the normalization constraint $\int_\R\rho_g=1$. 

\begin{figure}[h!btp] 
    \includegraphics[scale=0.6]{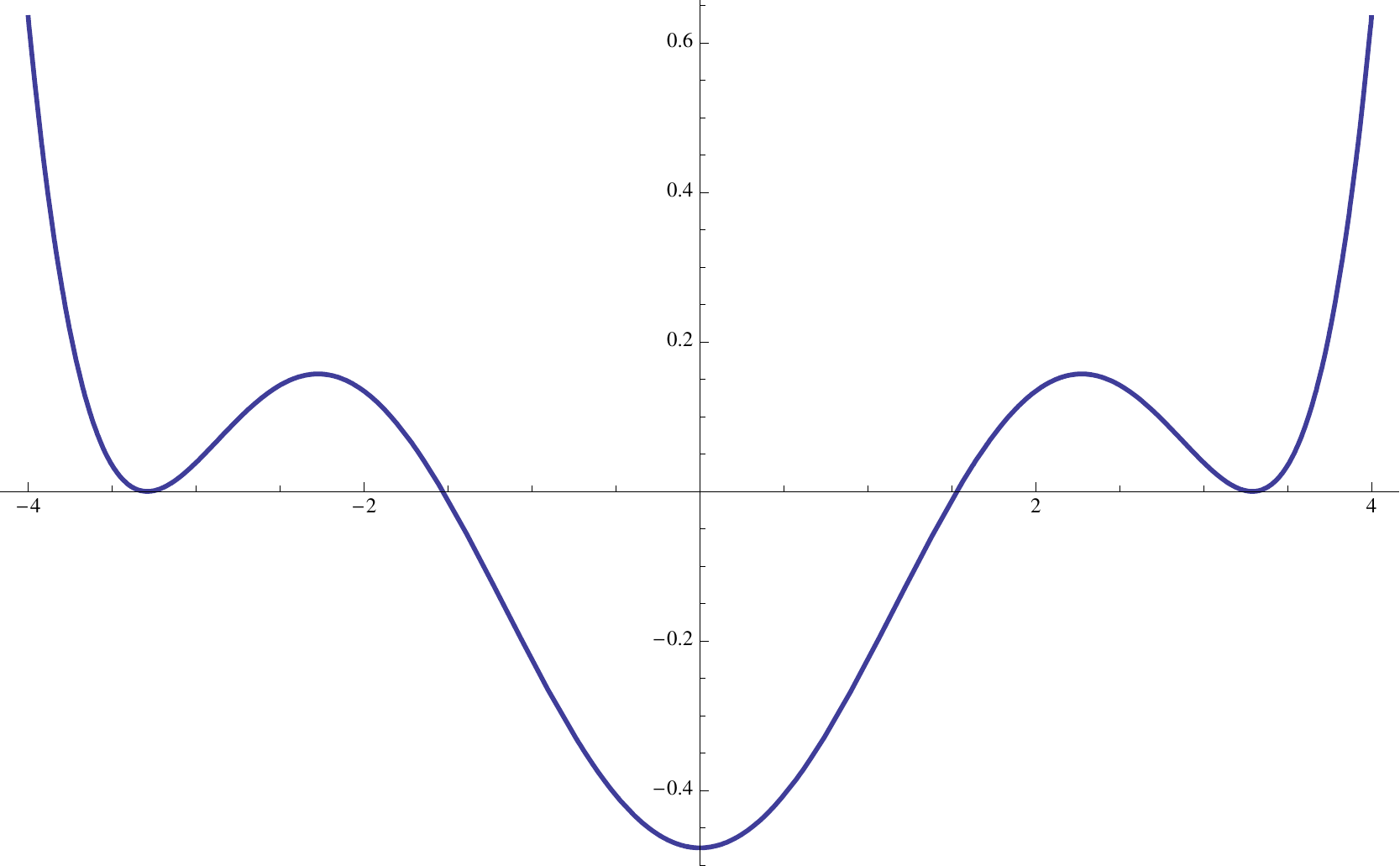}
     \caption{Polynomial $P_g(x)$ for $g=-1/48$ and $\beta=1$. }\label{polynomialP}
\end{figure}

Physical arguments lead us to the value of $J_g$. 
The six real roots of $P$ can not be all of multiplicity one otherwise the measure $\mu$ would have a disconnected 
support (union of three disjoints intervals), which is counter intuitive. Moreover the polynomial $P$ is even 
so that there are two roots with multiplicity two and two roots with multiplicity one. 
The zeroes with multiplicity two are symmetric and exterior while the zeroes with multiplicity one stand inside 
(see Fig. \ref{polynomialP}). 
We can compute the constant $J_g$ from the condition $P_g(\xi_g)=0$. We obtain 
\begin{align*}
J_g= -\frac{1}{\beta} \left(\left(2 g \xi_g^3 +\frac{\xi_g}{2}\right)^2- 2\beta g \xi_g^2  \right)\,. 
\end{align*}
The polynomial $P:=P_g$ is characterized completely and we can factorize it as
\begin{align*}
P_g(z)= 4 g^2 \, (z-\xi_g)^2 \, (z+ \xi_g)^2 \, (z^2 - \gamma^2)
\end{align*}
where 
\begin{equation}\label{def-gamma}
\gamma^2= - \frac{1}{6 g} \left( 1 - \sqrt{1+ 24 \beta g} \right)\,. 
\end{equation}
The Stieltjes transform $G_g$ is also determined: for $z\in\H$ near the real axis, we obtain  the 
following (explicit) expression for $G_g(z)$ after further (elementary) computations, 
\begin{align}\label{final-formula-Gg}
G_g(z) = \frac{2}{\beta} \left[-\frac{z}{2} - 2 g z^3 + 
\left(2 g z^2 + \frac{1}{6} \sqrt{1+ 24 \beta g} +\frac{1}{3}\right) 
 \sqrt{z^2-\gamma^2}   \right]\,. 
\end{align}

From the Stieltjes inversion formula, we can recover the probability density $\rho_g$, which  
is supported on the compact interval $[-\gamma;\gamma]$ where $\gamma>0$ is given in \eqref{def-gamma}. 
For $\lambda\in [-\gamma;\gamma]$, we obtain  
\begin{align}\label{density-super-quartic}
\rho_g(\lambda) = \frac{2}{\beta\pi} \left(2 g \lambda^2 + \frac{1}{6} \sqrt{1+ 24 \beta g} +\frac{1}{3}\right) 
 \sqrt{\gamma^2 - \lambda^2}\,. 
\end{align}
Reciprocally, we can check with an elementary integration that $\rho_g$ is 
 a probability measure and (using the residue Theorem) 
that its Stieltjes transform is indeed the analytic function $G_g$ 
characterized for $z$ near the real line in \eqref{final-formula-Gg}.

At the critical value $g_c=-1/(24\beta)$, we find 
\begin{align*}
P_{g_c}(z) = 4 {g_c}^2 (z- 2\sqrt{\beta})^3 (z+2\sqrt{\beta})^3\,, 
\end{align*} 
and
\begin{align}\label{density-critical-quartic}
\rho_{g_c}(\lambda) = \frac{1}{6 \pi \beta^2}  (4\beta-\lambda^2)^{3/2} \,. 
\end{align} 
For $\beta=2$, we re obtain the solution found in \cite[see their Figure 1]{bipz}. 

\begin{figure}[h!btp] 
    \includegraphics[scale=0.8]{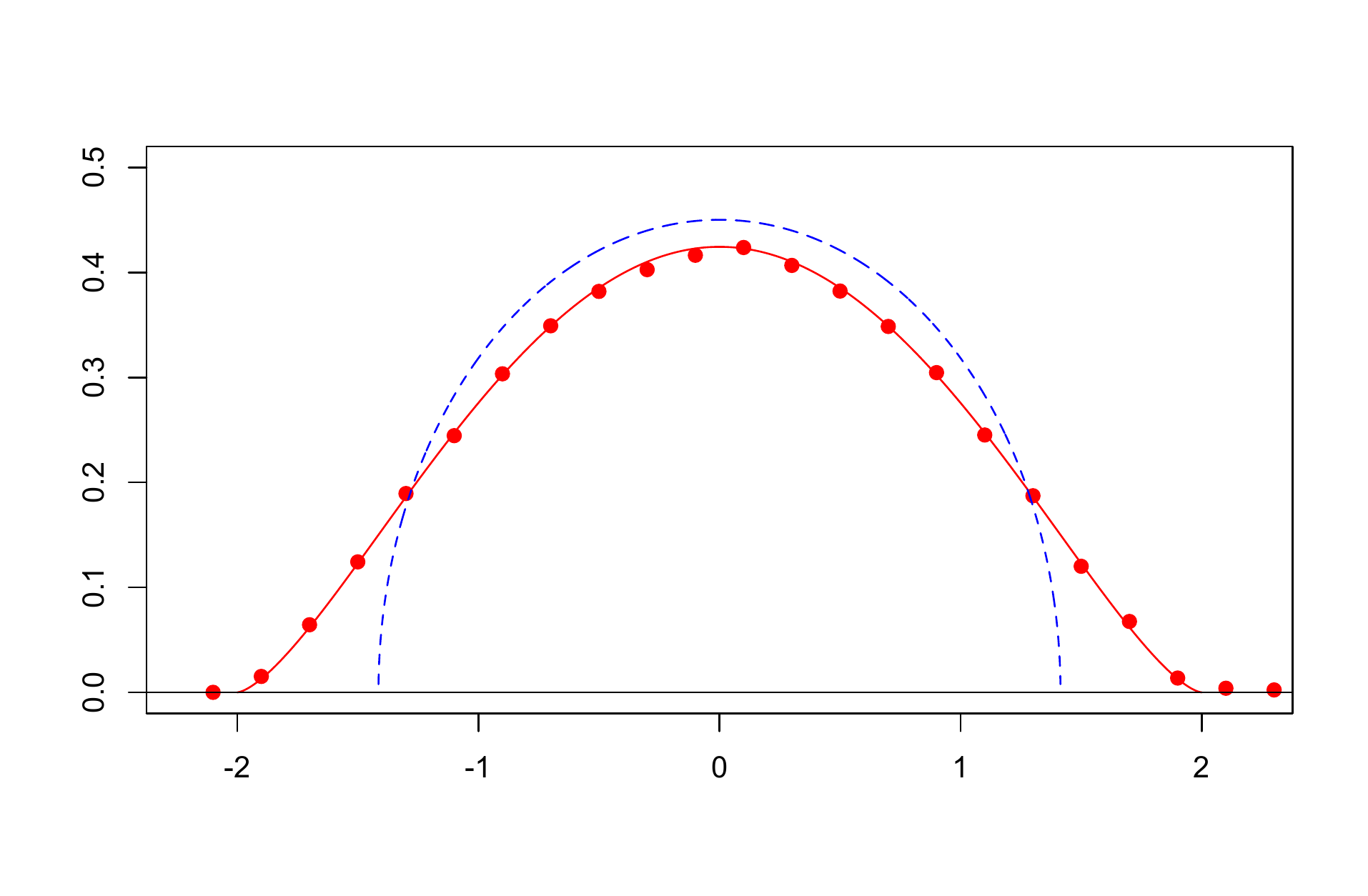}
     \caption{(Color online). Spectral density $\rho_g(x)$ for $g=g_c=-1/24$ and $\beta=1$ together with the semi-circle obtained for $g=0$. 
     The red points represent the empirical density obtained from our simulated samples in the critical case.} \label{quartic-density-fig}
\end{figure}

For $g\in [-\frac{1}{24\beta};0)$, the probability density functions $\rho_g$,  
given in \eqref{density-critical-quartic} for $g=g_c$ 
and \eqref{density-super-quartic} for $g> g_c$, are the limiting spectral density of the matrix $H(t)$ which follows the Langevin equation 
\eqref{langevin-quartic} in the stationary state as $N\to \infty$.  

For $g< g_c$, the situation is more complicated than that. One can still try to find the unique constant $J_g$ 
such that there exists $G_g$ with the two required properties and 
compute the imaginary part $\Im \, G_g(\lambda)/\pi$ for $\lambda\in \R$. 
But, in this case, one can check (at least numerically) that the function $\rho_g(\lambda):= \Im \, G_g(\lambda)/\pi$ 
is not a probability density. Its integral over $\R$ is in fact strictly smaller than $1$.
This simply means that there is not a constant $J\in \C$ such that 
the Stieltjes transform of some (probability) measure satisfies the 
quadratic equation \eqref{Stieltjes-quartic-equilibrium} for all $z\in \H$.
 
There is a rather clear physical interpretation to the non existence of such a constant $J$. 
For the cubic potential, the constant $J_a\in \H$ 
was related to the flux of probability in the system. 
The positive imaginary of $J_a$ was measuring the amount of probability density 
flowing  from the right to the left of $x\in\R$ per unit of time in the stationary state. In the present case of the Quartic 
potential with two sided exits, 
there can not be such a stationary flux of charges as the particles can exit 
either to the right by exploding in $+\infty$ or to the left in $-\infty$.  

It would actually be very interesting to  compute the 
limiting density $\rho_g$ of eigenvalues in the stationary state for values of $g<g_c$, by analyzing directly
\eqref{equilibrium-mu-quartic}.   
We would expect some diverging terms coming in at the origin, where there is a birth process
of new particles. 

\section{Conclusion and opening}\label{conclusion}

We now conclude this work with some open questions and further comments on possible extensions. 

We note that the case of a cubic potential was also considered at the end of \cite{bipz}. 
The authors study divergent 
matrix integrals of the form $\int \exp(-N W_g(M)) dM$ where $W_g(x)= x^2/2 + g x^3/3, g\in \R$ 
viewed as power series 
and they recover the
limiting spectral density when the confinement is 
strong enough so that the eigenvalues are localized in a compact support. 
The ideas developed in this paper permit one to define invariant ensembles associated to the potential $W_g$ for $g>0$. 
The stationary spectral densities can also be obtained in explicit forms, 
as was done in \cite{bipz} for $g$ smaller than a critical value $g_c$. 
The novel contribution of our work is that we went further by 
computing the spectral density in the regime $g>g_c$, where it has unbounded 
support with a stationary macroscopic current of particles across the system.  

More generally, the construction explained in section \ref{cubic-hermitian-process} permits one to define invariant ensembles 
of random matrices in general polynomial potentials $V(x)$ of arbitrary odd degree $k$, 
with possibly multiple wells. The Stieltjes transform method to study the large $N$ limit of the spectrum 
does not adapt straightforwardly though 
because unknown moments of the limiting spectral densities 
arise in the evolution equation satisfied by $G$ if ${\rm deg}(V) \geq 5$.
Those unknown moments have to be determined from the constraints and analyticity conditions but 
we expect the explosions to occur sufficiently fast 
for the tails of the spectral density to be light enough for those (few) moments to be finite. 
We also conjecture that such potentials lead to stationary spectral densities with support on several 
disconnected intervals, in the semi-confining regimes.  
A sharp phase transition should be observed as the confinement 
gets weaker, leading to spectral densities with unbounded support and
 associated to macroscopic stationary currents of particles in the system. 
In the semi-confining cases, the eigenvalues would lie in the wells near the local minimums of the potential $V$, leading 
to a more complex deformation towards the unbounded probabilities found in the fully non-confining regimes. 
We are currently working on this problem.  

It seems that other types of phase transitions for the spectral density at the critical value could be observed. 
Indeed, one can obtain, through an appropriate choice of the (polynomial) potential $V$ \cite{neuberger}, 
a limiting spectral density 
which vanishes at the end point of the spectrum $\lambda_c$ as $(\lambda-\lambda_c)^{k-1/2}$ for any integer $k \geq 1$. 
The potential $V$ is then called $k$ multi-critical. A general potential usually leads to criticality $k=1$, as 
in the case of the quadratic 
(Gaussian) potential. 
An interesting extension of our work would be to exhibit such (polynomial) potentials with odd degree which 
would cover general criticality parameters $k\geq 3$ and lead to more complex phase transitions 
for the spectral density. 
In our work, the potential $V_a$ at the critical value $a=a^*$ is $2$ multi-critical, 
leading to a density vanishing at the lower edge as 
$(\lambda- \lambda_c )^{3/2}$ (with $\lambda_c=-\frac{1}{2}\beta^{1/3}$). 

Other interesting questions would be to study the top eigenvalues statistics in the subcritical case $a< a^*$ when the 
density of eigenvalues has full support with tails $\rho_a(x) \sim C_a /x^2$. 
In particular, the fluctuations of the top eigenvalue around its asymptotic value of order $N$ are of interest in this fully non-confining regime. 
One may wonder whether those statistics are related to those of the top eigenvalues of heavy tailed Wigner matrices, or 
if they are of a different nature. In particular, the crossover for those statistics at the critical value $a=a^*$ 
is interesting. 
The eigenvalues statistics at the edge of the spectrum in the critical $a=a^*$ and super critical $a>a^*$ cases 
have been extensively investigated since $1990$. 
In $1991$, Bowick and Br\'ezin \cite{bb} have computed the spectral density
at the edge of the spectrum for general multi-critical potentials, in the scaling region of width $N^{-\nu}$ 
with $\nu=\frac{2k}{2k+1}$, 
where $k$ is the criticality (see also \cite{forrester2}). 
Their results are in fact universal in the sense that the scaling shape of the density at the edges 
does not depend on the particular choice of the potential $V$ of given criticality $k$.  
In particular, 
 for $k=1$ their work is pioneer and provides the first description of  the Tracy-Widom region.
In $1993$, Tracy and Widom \cite{tracy-widom1,tracy-widom2} proved 
the weak convergence of the top eigenvalue for the classical Gaussian ensembles 
(which correspond to the quadratic potential) to the so called Tracy-Widom distributions. 
A more precise description was later provided in \cite{virag} where the authors prove that the joint convergence of the 
top eigenvalues of the Gaussian ensembles to those of the stochastic Airy operator.  
The multi-critical cases $k\geq 2$ were later considered in \cite{claeys} and the results of \cite{virag} 
are conjectured to extend for those potentials as well (see the end of \cite{rider}). 
 
We conclude this paper with a last open question on the fluctuations of the random flux of particles 
around its asymptotic value given in Eq. \eqref{leading-flux}, 
in the stationary state and in the large $N$ limit. It would be very interesting to relate the order of the fluctuations, the shape 
of the fluctuations in the central regime or the large deviations regimes to 
other models of interacting particles with random flux.

\section{Proof of Theorem \ref{main} }\label{proof}
\subsection{Evolution equation of $G_N(z,t)$}\label{evolution.GN}

Recall the notations $\tau_1 <\tau_2< \cdots < \tau_k\cdots , k \in \N$ for the successive explosions times 
of the diffusive matrix process $H$. 

Following \cite{rogers,jp-alice,satya}, 
we look for an evolution equation for the Stieltjes transform $G_N(z,t)$. By It\^o's formula, we have for any $z\in \H$ 
and for any $t\in (\tau_k;\tau_{k+1})$, 
\begin{align}\label{burger.brut}
dG_N(z,t) = - \frac{a}{N} \sum_{i=1}^N \frac{dt}{(\lambda_i-z)^2} + \frac{1}{N} 
\sum_{i=1}^N \frac{\lambda_i^2}{(\lambda_i-z)^2} dt 
- \frac{\beta}{2N^2} \sum_{i\neq j} \frac{1}{(\lambda_i-z)^2} \frac{dt}{\lambda_i-\lambda_j} \\ + 
\frac{1}{2N^2} \sum_{i=1}^N \frac{2}{(\lambda_i-z)^3} dt  - \frac{1}{N^{3/2}} \sum_{i=1}^N \frac{dB_i}{(\lambda_i-z)^2} \,.  
\end{align}

At the explosion times $\tau_k,k\in \N$, the eigenvalue which explodes in $-\infty$, is immediately restarted at $+\infty$, 
so that the Stieltjes transform $G_N(z,t)$ of the empirical measure
is continuous in $\tau_k$. The evolution equation \eqref{burger.brut} therefore holds for all $t\geq 0$.  

We want to rewrite this equation in a more compact way as a function of $G_N$. 
Using symmetry and anti-symmetry properties, we have  
\begin{align*}
 &\sum_{i\neq j} \frac{1}{(\lambda_i-z)^2} \frac{1}{\lambda_i-\lambda_j} = \frac{1}{2} 
  \sum_{i\neq j}  \frac{1}{\lambda_i-\lambda_j} \left( \frac{1}{(\lambda_i-z)^2}- \frac{1}{(\lambda_j-z)^2}\right)\\
&= \frac{1}{2}  \sum_{i\neq j}  \frac{1}{\lambda_i-\lambda_j} \frac{(\lambda_j-\lambda_i)^2 + 2 (\lambda_j-\lambda_i)(\lambda_i-z)}{(\lambda_i-z)^2(\lambda_j-z)^2}=- \sum_{i\neq j} \frac{1}{(\lambda_i-z)(\lambda_j-z)^2}\,. 
\end{align*}
The last expression conveniently rewrites in terms of $G_N$ as 
\begin{align*}
\sum_{i\neq j} \frac{1}{(\lambda_i-z)(\lambda_j-z)^2} &= \sum_{i,j=1}^N \frac{1}{(\lambda_i-z)(\lambda_j-z)^2}
- \sum_{i=1}^N \frac{1}{(\lambda_i-z)^3} \\
&=N^2 \, G_N \, \partial_z G_N -  \frac{N}{2} \partial_z^2 G_N\,. 
\end{align*}
From \eqref{burger.brut} and the previous computations, we have, for any $z\in \H$ and $t\geq 0$,
\begin{align}\label{burgers.net}
 G_N(z,t) - G_N(z,0) &=  \int_0^t  \partial_z \left[ \frac{\beta}{4}  \, G_N^2(z,s) + (z^2-a) \, G_N(z,s) + z +  
\frac{1}{2N}(1- \frac{\beta}{2}) \partial_z G_N(z,s) \right]   ds  \notag \\  
&- \frac{1}{N^{3/2}} \sum_{i=1}^N \int_0^t \frac{dB_i(s)}{(\lambda_i(s)-z)^2}\,. 
\end{align}
We aim at taking $N\to \infty$. 
For this, we need to prove that the sequence of probability measures process $\{ (\mu_t^N)_{\{0 \leq t \leq T\}}, N \in \N\} $ is
almost surely pre-compact in the space $\mathcal{C}([0,T],\mathcal{M}_{\le 1}(\R))$ where $\mathcal{M}_{\le 1}(\R)$ is 
the space of Borel measures $\nu$ on $\R$ with total mass $\nu(\R) \leq 1$ equipped with its weak-$\star$ topology.

It is convenient to work with the weak-$\star$ topology on $\mathcal{M}_{\le 1}(\R)$ (where a sequence $\mu_n$ converges to $\mu$ iff $\int f d \mu_n$ converges to $\int f d \mu$ for all $f\in C_c(\R)$, the space of continuous and compactly supported functions on $\R$) as this topology is metrizable on this bounded set and makes $\mathcal{M}_{\le 1}(\R)$ compact\footnote{It is closed and sequentially compact}. Therefore, the topology on $\mathcal{C}([0,T],\mathcal{M}_{\le 1}(\R))$ we consider is simply the topology of uniform convergence for this metric. Note that the weak-$\star$ topology on the bounded set $\mathcal{M}_{\le 1}(\R))$ is equivalent to the vague topology\footnote{But of course, as we are working with sub-probability measures, the weak-$\star$ topology is weaker than the usual topology of weak-convergence where the limit should hold for all continuous bounded functions.} (where $\mu_n$ converges to $\mu$ iff $\int f d \mu_n$ converges to $\int f d \mu$ for all $f\in C_0(\R)$, the continuous functions on $\R$ which tends to $0$ when $x \to \pm \infty$).
 
Let us emphasize that, in contrast with the usual case handled in \cite{agz,rogers} 
where the authors prove a dynamical version of Wigner's Theorem, we do not 
need to prove pre-compactness in the smaller space $\mathcal{C}([0,T],\mathcal{P}(\R))$ as we will prove that the continuous limiting process (along subsequences) necessarily takes values in the space of probability measures.

This almost sure pre-compactness in $\mathcal{C}([0,T],\mathcal{M}_{\le 1}(\R))$ is proved in Lemma \ref{precompactness}. 
From any subsequence of $\N$, we can now extract a sub-sub-sequence $(N_k)_{k \geq 0 }$ such that we have the following 
pointwise convergence in $\mathcal{C}([0,T],\mathcal{M}_{\le 1}(\R))$, 
\begin{align*}
(\mu_t^{N_k})_{0\leq t \leq T} \Rightarrow (\mu_t)_{0 \leq t \leq T}\,,
\end{align*}
where $(\mu_t)_{0 \leq t \leq T} \in \mathcal{C}([0,T],\mathcal{M}_{\le 1}(\R))$. 

To any limit point $ (\mu_t)_{0 \leq t \leq T}$, we associate its Stieltjes transform process $(G(\cdot,t))_{0 \leq t \leq T}$
such $G(\cdot,t)$ is the Stieltjes transform of the measure $\mu_t$ for any $t\geq 0$. 

We note that, for $z\in \H$ fixed, the function $f(x)=1/(x-z)$ 
and its derivative $f'(x)= -1/(x-z)^2$ are continuous and tends to $0$ when $x \to \pm \infty$ and therefore belongs to $C_0(\R)$. We deduce the following convergences 
(along the subsequence $(N_k)$) for any $t\in  [0,T]$, 
 \begin{align*}
G_N(z,t) \longrightarrow G(z,t)\,, \quad  \partial_z G_N(z,t) \longrightarrow \partial_z G(z,t)\,. 
\end{align*}
Note in addition that the martingale term 
\begin{align*}
M_t^N = \frac{1}{N^{3/2}} \sum_{i=1}^N \int_0^t \frac{dB_i(s)}{(\lambda_i(s)-z)^2}
\end{align*}
which appears in \eqref{burgers.net} has a quadratic variation $\langle M_{\cdot}^N \rangle_t$ smaller then 
$t/(N^2 |\Im z|^{4})$.  By the Burkholder-Davis-Gundy inequality \cite[Theorem H.8]{agz} 
and the Chebyshev's inequality, we get that, for a universal constant $C>0$,
\begin{align}\label{martingale-control}
\P[\sup_{0\leq t\leq T} |M_t^N| \geq \varepsilon ] \leq \frac{C \, T}{\varepsilon^2 N^2|\Im z|^{4} }\,. 
\end{align}
By the Borel-Cantelli Lemma, for any $z\in \H$, we have $\sup_{0\leq t\leq T} |M_t^N|\to 0$ as $N\to +\infty$ almost surely. 
 
Eventually,  from \eqref{burgers.net}, we obtain the following equation 
satisfied by the Stieltjes transform process $(G(\cdot,t))_{0 \leq t\leq T}$ of 
any limit point $ (\mu_t)_{0 \leq t\leq T}$ of the pre-compact family 
$((\mu_t^N)_{0 \leq t\leq T}), N \in \N$ (see Lemma \ref{precompactness-lemma}) ,  
\begin{align}\label{burgers-eq2}
G(z,t) = G(z,0) + \int_0^t \partial_z  \left[   \frac{\beta}{4}  \, G(z,s)^2 + (z^2-a) \, G(z,s) + z  \right] \, ds \,, 
\end{align}
which holds for all $t\in [0;T], z\in \H$. 

It remains to prove that the measure $\mu_t$ is indeed a probability measure for any $t\geq 0$. 
Let $t \geq 0$, $\varepsilon >0$ and $z= \mathbf{i} /\varepsilon^2\in \H$.  
From \eqref{burgers-eq2}, we easily deduce that for all $\varepsilon >0$,  
\begin{align}\label{eq.G-eps}
&G(\frac{\mathbf{i}}{\varepsilon^2},t+\varepsilon) - G(\frac{\mathbf{i}}{\varepsilon^2},t)\\ & = 
\int_{t}^{t+\varepsilon} ds \left[\frac{\beta}{2} G(\frac{\mathbf{i}}{\varepsilon^2},s) \partial_z G(\frac{\mathbf{i}}{\varepsilon^2},s) - 
(\frac{1}{\varepsilon^4}+a) \partial_z G(\frac{\mathbf{i}}{\varepsilon^2},s) +  \frac{2\,\mathbf{i}}{\varepsilon^2}  G(\frac{\mathbf{i}}{\varepsilon^2},s) + 1 \right] \,.
\end{align}
Recalling $G(z,t)=\int_\R \frac{\mu_t(dx)}{x-z}$, we can check that the left hand side of \eqref{eq.G-eps} 
is of order $\varepsilon^2$ (or even smaller) as $\varepsilon \to 0$ 
while the right hand side is equivalent to $\varepsilon(1-\mu_t(\R))$ (using in addition the continuity of the function $t \mapsto \mu_t(\R)$). Therefore we have $\mu_t(\R)=1$ for every $t \geq 0$.

Finally, the convergence of $(\mu_t^N)_{0 \leq t \leq T}$ 
to $(\mu_t)_{0 \leq t \leq T}$ follows by uniqueness of the solution of \eqref{burgers-eq2} (proved in Lemma \ref{unicity-lemma}). 
Theorem \ref{main} is proved.

\subsection{Pre-compactness of the family $((\mu_t^N)_{t\geq 0}), N \in \N$}  \label{precompactness}

Denote by $\mathcal{M}_{\le 1}(\R)$ the space of Borel measures on $\R$ with total mass smaller or equal to $1$ equipped with its weak-$\star$ topology and 
recall that $\mathcal{P}(\R)$ denotes the space of probability measures on $\R$. 

 \begin{lemma}\label{precompactness-lemma}
 The family $\{ (\mu_t^N)_{\{0 \leq t \leq T \}}, N \in \N\} $ of continuous process in  $\mathcal{P}(\R)$ 
 is almost surely pre-compact in the space $\mathcal{C}([0,T],\mathcal{M}_{\le 1}(\R))$. 
 \end{lemma}
 
{ \it Proof. }
We first describe a family of compact subsets of $\mathcal{C}([0,T],\mathcal{M}_{\le 1}(\R))$. 
Let $(f_i)_{i\geq 0}$ be a sequence of bounded continuous functions dense in the space $\mathcal{C}_0(\R)$ of continuous 
function on $\R$ which converge to $0$ in $\pm \infty$ and let $(C_i)_{i \geq 0}$ be a family of compacts subsets of 
 the space $\mathcal{C}([0,T],\R)$ of continuous functions from $[0,T]\to \R$. Then, adapting the proof of Lemma 4.3.13 in \cite{agz}, we can prove that the set  
 \begin{align*}
\mathcal{K} := \bigcap_{i \geq 0} \{ t \to \int _\R f_i \, d\mu_t \in C_i \}  
\end{align*}
is a compact subset of $\mathcal{C}([0,T],\mathcal{M}_{\le 1}(\R))$. This proof is straightforward showing that $\mathcal{K}$ is closed and sequentially compact\footnote{The space $\mathcal{C}([0,T],\mathcal{M}_{\le 1}(\R))$ is metrizable.} with a diagonal extraction, noting 
in addition that $\mathcal{M}_{\le 1}(\R)$ is compact for the topology of weak-$\star$ convergence.  

We denote by $\mathcal{F}$ the set of twice differentiable functions $f:\R\to \R$, 
such that $||f||_\infty < \infty, ||f'||_\infty < \infty$ and $||f''||_\infty < \infty$
and such that, in addition, $||x^2 f'(x)||_\infty < \infty$. 
Note that $\mathcal{F}$ is dense in $\mathcal{C}_0(\R)$. 
We need an estimate on the H\"older norm of the function 
$t\to \int_\R f(x) \mu_t^N(dx)$ for any $f\in \mathcal{F}$.

Applying It\^o's formula, we get for any $s<t$,
\begin{align*}
 \int_\R f(\lambda) \mu_t^N(d\lambda) -&  \int_\R f(\lambda) \mu^N_s(d\lambda) = 
 \int_s^t \int_\R (a- \lambda^2) f'(\lambda) \mu^N_u(d\lambda) du +  \frac{\beta}{2 N^2} \int_s^t  \sum_{i\neq j} \frac{f'(\lambda_i(u))}{\lambda_i(u)-\lambda_j(u)} \, du \\
 &+\frac{1}{N^{3/2}} \sum_{i=1}^N f'(\lambda_i(t)) dB_i(t) + \frac{1}{2N} \int_s^t  \int_\R f''(\lambda) \mu^N_u(d\lambda) \, du \,. 
 \end{align*}
Using symmetry, we have 
\begin{align*}
 \sum_{i\neq j} \frac{f'(\lambda_i)}{\lambda_i-\lambda_j} = \frac{1}{2} \sum_{i\neq j}  \frac{f'(\lambda_i)-f'(\lambda_j)}{\lambda_i-\lambda_j}\,. 
\end{align*}
Therefore, 
\begin{align*}
\frac{\beta}{2 N^2}\left | \int_s^t  \sum_{i\neq j} \frac{f'(\lambda_i(u))}{\lambda_i(u)-\lambda_j(u)} \, du \right| 
\leq \frac{\beta}{4 N^2} N(N-1) ||f''||_\infty (t-s) \leq \frac{\beta}{4} (t-s) ||f''||_\infty \,. 
\end{align*}
Denote by $C_f := \max\{||f||_\infty, \;||f'||_\infty,\; ||f''||_\infty,\;||x^2\,f'(x)||_\infty\} < \infty$.
Gathering the above estimates, it follows that, such that for any $s,t\in [0,T]^2$,  
\begin{align*}
| \int_\R f(\lambda) \mu_t^N(d\lambda) -  \int_\R f(\lambda) \mu^N_s(d\lambda)  | \leq 
\left(a + \frac{\beta}{4} +\frac{3}{2}\right) \,C_f\; |t-s| + |M_f^N(t)- M_f^N(s)| 
\end{align*}
where $(M_f^N(t))_{t\geq 0}$ is the martingale process defined as
\begin{align*}
M_f^N(t) = \frac{1}{N^{3/2}} \sum_{i=1}^N f'(\lambda_i(t)) \, dB_i(t) \,. 
\end{align*}
With the same proof as in \cite[Proof of Lemma 4.3.14, page 265]{agz}, we can prove that there exists a 
constant  $C>0$ which depends only on $T$ such that, 
for any $\delta>0, M>0$, 
\begin{align}\label{continuity-martingale}
\P\left[\sup_{|t-s| \leq \delta \atop 0\leq s,t \leq T} |M_f^N(t)-M_f^N(s) | \geq M \delta ^{1/8}\right]  
\leq \frac{C \delta^{1/2}}{N^2 M^4} ||f'||_\infty^2\,.
\end{align}
Using \eqref{continuity-martingale}, we deduce that for any $\delta:= \delta(a,\beta,C_f)>0$ small enough and $M>0$, we have 
\begin{align}\label{continuity-density}
\P\left[\sup_{|t-s| \leq \delta \atop 0\leq s,t \leq T} |\int_\R f(\lambda) \mu_t^N(d\lambda) -  \int_\R f(\lambda) \mu^N_s(d\lambda)| \ge
 (M+1) \delta ^{1/8} \right]  \leq \frac{C \delta^{1/2}}{N^2 M^4}||f'||_\infty^2\,. 
\end{align}
Recall that, by the Arzela-Ascoli Theorem, sets of the form 
\begin{align*}
C:= \bigcap_{k \in \N} \{ g \in \mathcal{C}([0;T], \R) :
 \sup_{|t-s| \leq \delta_k \atop 0\leq s,t\leq T} |g(t)-g(s)| \leq \epsilon_k , \sup_{0\leq t \leq T} |g(t)| \leq 1/\alpha_k \}\,. 
\end{align*}
where $(\delta_k)_{k},(\epsilon_k)_k,(\alpha_k)_k$ are sequences of positive real numbers going to zero as 
$k$ goes to infinity, are compact in $\mathcal{C}([0;T], \R)$.  

For $f\in \mathcal{F}$ and $\varepsilon>0$, we consider the subset 
of $\mathcal{C}([0;T], \mathcal{M}_{\le 1}(\R))$, defined by 
\begin{align*}
C_T(f,\varepsilon):= \bigcap_{k \in \N} \{ (\nu_t) 
\in\mathcal{C}([0;T], \mathcal{M}_{\le 1}(\R)) :  
\sup_{|t-s| \leq k^{-4} \atop 0\leq s,t\leq T} |\int_\R f(x) \nu_t(dx) -\int_\R f(x) \nu_s(dx) | 
\leq \frac{1}{\varepsilon \sqrt{k}} \}\,. 
\end{align*}
Then, using \eqref{continuity-density}, we have 
\begin{align*}
\P\left[(\mu_t^N)_{t\in [0;T]} \in C_T(f,\varepsilon)^c \right] \leq \frac{C\varepsilon^4}{N^4}\,. 
\end{align*}
We now pick a dense family $(f_i)_{i\in \N}$ in $\mathcal{C}_0(\R)$ of functions $f_i\in \mathcal{F}$ 
and setting $\varepsilon_i= 1/i$, we define the subset 
\begin{align*}
\mathcal{K} := \bigcap_{i \in \N} C_T(f_i,\varepsilon_i) \subset 
\mathcal{C}([0;T], \mathcal{M}_{\le 1}(\R))\,.
\end{align*}
By the Borel-Cantelli Lemma, we get 
\begin{align*}
\P[\bigcup_{N_0\in \N} \bigcap_{N\geq N_0} (\mu_t^N)_{t\in [0;T]} \in \mathcal{K} ] = 1\,,
\end{align*}
and the Lemma follows since $\mathcal{K}$ is a compact set of $\mathcal{C}([0;T], \mathcal{M}_{\le 1}(\R))$. 

\qed

\subsection{Proof of uniqueness of solutions of \eqref{burgers-eq}}  
 \begin{lemma}\label{unicity-lemma}Let $\mu\in \mathcal{P}(\R)$. 
There exists a unique (deterministic) process 
$(G(\cdot,t))_{t\geq 0}$ in the space of analytic function $\H\to \H$ 
 which enjoys the two following properties: 
\begin{itemize}
\item for any $t\geq 0$, $G(\cdot,t)$ is the Stieltjes transform of a real probability measure;  
\item $G$ is a strong solution of the holomorphic partial differential equation 
on $\H\times \R_+$,
\begin{align}\label{burgers-eq3}
\partial_t G(z,t) = \partial_z  \left[   \frac{\beta}{4}  \, G(z,t)^2 + (z^2-a) \, G(z,t) + z  \right]\,, 
\end{align}
with initial condition $G(z,0)= \int_\R \frac{\mu(dx)}{x-z}$. 
\end{itemize} 
  \end{lemma}
  
{\it Proof.} We set 
 \begin{equation}\label{relation-G-H}
 H(z,t):=G(z,t) +\frac{2}{\beta} (z^2-a)\,. 
 \end{equation}
 It is straightforward to check that the analytic function $H(\cdot,t)$ satisfies the evolution equation for $t\geq 0$, 
 \begin{align}\label{pde-H}
\partial_t H(z,t) = \frac{\beta}{2} H(z,t) \partial_z H(z,t) - \frac{4}{\beta} z (z^2-a) +1 \,. 
\end{align} 
It suffices to prove that there is at most one unique strong solution $H:\H \times \R_+ \to \C$ of the partial differential equation (pde) 
\eqref{pde-H}
 such that in addition $H(\cdot,t):\H\to\H$ is analytic for any $t\geq 0$.  Following \cite{agz,rogers}, we use the characteristic method. 
For $z_0\in \H$, we consider the following Cauchy problem 
\begin{align}\label{cauchy-pb}
z'(t) = - \frac{\beta}{2} H(z(t),t)\,, \quad z(0) = z_0. 
\end{align} 
For $\eta >0$ and any $z_1,z_2\in \H^\eta:=\{z\in \C: \Im z > \eta\} $, we have 
for any $t\geq 0$,
\begin{align*}
|G(z_1,t)-G(z_2,t)| \leq \int_\R \mu_t(dx) |\frac{1}{x-z_1}-\frac{1}{x-z_2}| \\
\leq |z_1-z_2| \int_\R \frac{\mu_t(dx)}{|x-z_1||x-z_2|} \leq \frac{|z_1-z_2|}{\eta^2}\,.  
\end{align*}
We deduce that the continuous function $(z,t)\in \H^\eta\times \R_+ \to H(z,t)$ is {\it locally} in $z$ (globally in $t$) Lipschitz 
on $\H^\eta\times \R_+$ with respect to the first variable $z$. More precisely, we mean that, for any compact set $K\subset \H^\eta$, 
there exists a constant $M>0$ such that for any $t\geq 0$ and $z_1,z_2\in K$, we have 
\begin{align}\label{cond-Lip-H}
| H(z_1,t) - H(z_2,t) | \leq M \, |z_1- z_2|\,. 
\end{align}
By the Cauchy Lipschitz theorem, for any $z_0\in \H^\eta$, there exists a unique solution $(z(t))_{0\leq t < \varepsilon}$
to the Cauchy problem \eqref{cauchy-pb}, defined up to a time $\varepsilon >0$ small enough such that $z(t) \in \H$ 
for all $t< \varepsilon$ 
(recall that $H(\cdot,t)$ is analytic on $\H$ so that the solution has to remain in this domain for the differential equation in \eqref{cauchy-pb}
to be well defined).

Let us now fix $\omega_0\in \H^\eta$. 
We can prove using the Lipchitz condition \eqref{cond-Lip-H} and \eqref{cauchy-pb}
that we can pick $r>0$ small enough (depending only on $\omega_0$ and on the constant $M$) such that {\it any} solution of the Cauchy problem 
starting from any $z_0$ in the open ball $B(\omega_0,r)$ with center $\omega_0$ and radius $r>0$ is defined up to a time $\varepsilon >0$ 
(which depends on $M$ and $r$ but which is \emph{independent} of $z_0\in B(\omega_0,r)$).

Differentiating  \eqref{cauchy-pb} with respect to $t$ and using also the pde \eqref{pde-H}, we obtain a new {\it explicit} 
Cauchy problem for the function $z$, 
\begin{align}\label{cauchy-pb2}
z''(t) = 2 z(t) (z^2(t)-a) - \frac{\beta}{2} = \frac{1}{2} P'(z(t)), \quad z(0)=z_0, \quad z'(0)  = - \frac{\beta}{2} H(z_0,0)\,,
\end{align} 
where $P'$ is the polynomial we have already met in \eqref{Pprime}.

We now regard the Cauchy problem \eqref{cauchy-pb2} as a holomorphic Cauchy problem on the domain 
$(t,z) \in \C \times \C$. We know from the fundamental theorem \cite{yakovenko} that there exists a unique solution to \eqref{cauchy-pb2}
defined for all $t\in \C$ and for any initial condition $z_0\in \H^\eta$. Moreover, $z(t)$ depends holomorphically on the initial conditions: there exists 
a holomorphic function $F_t$ such that, for any $t\in \C$, 
\begin{align*}
z_t = F_t(z_0)\,.  
\end{align*}
By uniqueness, for any $z_0\in B(\omega_0,r)$, the solution $(z(t))$ of the Cauchy problem \eqref{cauchy-pb2} coincides 
with the one of the Cauchy problem \eqref{cauchy-pb} on the small interval $[0,\varepsilon )$. 

We now fix $t\in [0,\varepsilon]$ and we consider the {\it open and simply connected domain} $\Omega:=F_t(B(\omega_0,r))$ 
of the ``targets". By definition of $\varepsilon$, we have $\Omega\subseteq \H$.  
Besides, by uniqueness, the solutions of the Cauchy problem \eqref{cauchy-pb} can not coincide in $\H$ at the given time $t$, 
and $F_t : B(\omega_0,r) \to F_t(B(\omega_0,r))$ 
is a conformal isomorphism (i.e. analytic and bijective). 

Finally, for any $t \in [0,\varepsilon]$ and $\omega\in \Omega$ fixed, we can find a unique $z_0\in B(\omega_0,r)$ such that the solution 
of \eqref{cauchy-pb} satisfies $z(t)= \omega$. 
Then the formula 
\begin{align*}
H(\omega,t) = - \frac{2}{\beta} \,  z'(t)\,,
\end{align*}
valid for any $\omega$ in the open and simply connected domain $\Omega$,  
characterizes uniquely the analytic function $H(\cdot,t) : \H \to \H$. 
This argument is true for any $t\leq \varepsilon$ so that $H(z,t)$ is uniquely characterized for all $z\in \H$ and $t\in [0;\varepsilon]$.  
 
 The same method (starting from the same initial condition $\omega_0$ for example so that the same $\varepsilon$ will work) permits us to extend this characterization on the intervals $[\varepsilon,2\varepsilon]$, $[2\varepsilon,3\varepsilon]$, ... 
 
 Uniqueness of $G$  is implied by the uniqueness of $H$ from \eqref{relation-G-H}. 
 
 \qed
 
 \subsection{Proof of Proposition \ref{prop-conv-time}}\label{proof-prop-conv-time}
We first prove the existence part of (1). 
It is easy to see that the solution $G$ of the evolution equation \eqref{burgers-eq} 
is Lipchitz in $t$ locally uniformly in $z$ i.e. for all compact subset $K\subseteq \H$, there exists a constant $M$
such that for all $t\geq 0$ and all $z\in K$, 
\begin{align*}
|G(z,t) - G(z,s)| \leq M \, |t-s|\,. 
\end{align*}
Therefore, using the Cauchy formula, we deduce that the holomorphic function $z \mapsto \partial_z[\beta G(z,s)^2/4 + (z^2 - a)G(z,s) + z]$ is Lipchitz in $t$ locally uniformly in $z$ as well.
The existence of the limit when $t \to \infty$ of $G(z,t)$ leads to the existence of the limit of the integral on the right hand side of \eqref{burgers-eq}. As the integrated function is uniformly continuous, we deduce that 
$G_a$ is a stationary solution of \eqref{burgers-eq}, i.e. 
\begin{align*}
\partial_z  \left[   \frac{\beta}{4}  \, G_a(z)^2 + (z^2-a) \, G_a(z) + z  \right] = 0\,,
\end{align*}
so that there exists a constant $J\in \C$ such that for all $z\in \H$, 
\begin{align}\label{eq-statio-G-infty}
\frac{\beta}{4}  \, G_a(z)^2 + (z^2-a) \, G_a(z) + z = J\,. 
\end{align} 
From \eqref{def-G-infty}, we already know that $G_a$ is bounded 
in the neighbourhood of $\infty$. Therefore, we deduce from \eqref{eq-statio-G-infty} that 
indeed $\lim_{y \to \infty} i y \,  G_a(i y) =1$. 
We then know from \cite[Theorem 1]{geronimo} that $G_a$ is the Stieltjes transform of a probability measure $\mu_\infty$
and that $\mu_t$ converges weakly to $\mu_\infty$. 
The proposition is proved.  
  
\appendix 

\section{Boltzmann weight of the Hermitian diffusion process $H$} 
Let us check that the probability distribution $P$ defined in \eqref{inv-ens} is 
a stationary measure of the stochastic differential system
\eqref{langevin-matrix}.  First notice that, if $M$ is a $N\times N$ real matrix, 
then the gradient of the function $M\to {\rm Tr}(V(M))\in \R$ with respect to the $N^2$ entries of 
the matrix $M$, is the function $M\to V'(M^\dagger)$. Thus, 
if $H$ is a Hermitian matrix, we simply have 
\begin{align}\label{nice.equality}
\nabla \, {\rm Tr}(V(H)) = V'(H)\,. 
\end{align}
 It remains to check that the probability distribution $P$ is the unique stationary solution of the {\it Fokker Planck}
 equation satisfied by the (stationary) transition probability of the diffusion process $(H(t))$,  
\begin{align}\label{fokker-planck}
\frac{\partial P}{\partial t} = 0 
= \frac{1}{2}\, \sum_{i,j=1} \frac{\partial}{\partial H_{ij}} \left[\left(\nabla {\textrm{Tr}}(V(H))\right)_{ij} P(H) \right] + \frac{1}{2N} 
\sum_{i,j=1} \frac{1+ \delta_{i=j}}{2} \frac{\partial^2}{\partial H_{ij}^2} P(H)\,. 
\end{align}  
The reader may actually check that the function $P$ as defined in \eqref{inv-ens} satisfies, for any Hermitian matrix $H$, the following conditions 
\begin{align}
\frac{\partial }{\partial H_{ii}} P(H) &=- N \, \left(\nabla {\textrm{Tr}}(V(H))\right)_{ii} P(H)\,,\notag \\
\frac{\partial }{\partial H_{ij}} P(H)= \frac{\partial }{\partial H_{ji}} P(H)&= -2 N \, \left(\nabla {\textrm{Tr}}(V(H))\right)_{ij} P(H) \, \quad {\mbox{if}} \quad i < j\,, \label{second.condition}
\end{align} 
under which \eqref{fokker-planck} trivially holds. The factor $2$ which appears in the second line \eqref{second.condition} 
is due to the symmetry of the matrix $H$.

\section{Stieltjes transform properties}   \label{stieltjes}
The Stieltjes transform is frequently used in random matrix theory for the study of empirical spectral densities in the large $N$ limit. 

A measure $\mu$ is characterized by its Stieltjes transform, which is an analytic function $G: \H \to \H$ ($\H$ denotes the open upper half-plane), defined as 
\begin{equation*}
G(z):=\int_\R \frac{\mu(dx)}{x-z}\,. 
\end{equation*} 
We have the following inversion formula valid for any measure $\mu$ on $\R$,
\begin{align}\label{inversion-formula}
\lim_{\varepsilon \downarrow 0}\int_x^y \Im \, G(\lambda+i \, \varepsilon) \, d\lambda= \pi\, \mu(x;y) + \frac{\pi}{2}(\mu(\{y\})-\mu(\{x\})) \,,  
\end{align}
where $\Im z$ denotes the imaginary part of $z\in \C$. 

When the Stieltjes transform $G(z)$ has a continuous extension to $\R \cup \H$, it is easy to check that $\mu$ admits a smooth density with respect to the Lebesgue measure.

If $\mu$ is a probability measure, its Stieltjes transform $G(iy)$ behaves as $-1/(iy)$ when $y$ goes to $+\infty$. 
Reciprocally, Akhiezer's theorem \cite[page 93]{akhiezer} states a useful criterium characterizing Stieltjes transforms of probability measure: $G$ is the Stieltjes transform of a probability measure iff $G$ is analytic on $\H$ with $G(\H) \subseteq \H$ and $G(iy) \sim -1/(iy)$ as $y \to + \infty$.

\end{document}